\documentclass[12pt,draftcls,onecolumn]{IEEEtran}

\usepackage{amsfonts}
\usepackage{amssymb}
\usepackage{amsmath}
\usepackage{mathptmx}
\usepackage{cancel}
\usepackage{pstricks}
\usepackage{psfrag}
\usepackage{mathrsfs}
\usepackage{graphicx}
\usepackage{pgf,tikz}
\usepackage{algorithm}
\usepackage{algpseudocode}
\usepackage{soul}
\usetikzlibrary{arrows}
\def\be{\begin{equation}}
\def\ee{\end{equation}}
\def\bea{\begin{eqnarray}}
\def\eea{\end{eqnarray}}
\def\beann{\begin{eqnarray*}}
\def\eeann{\end{eqnarray*}}

\def\ns{\hspace{-1mm}}
\newcommand{\real}{{\mathbb{R}}}

\def\spacingset#1{\def\baselinestretch{#1}\small\normalsize}
\spacingset{1.35}

\newtheorem{lemma}{Lemma}
\newtheorem{theorem}{Theorem}
\newtheorem{remark}{Remark}
\newtheorem{corollary}{Corollary}

\newtheorem{property}{Property}
\newtheorem{example}{Example}[section]

\def\be{\begin{equation}}
\def\ee{\end{equation}}
\def\bea{\begin{eqnarray}}
\def\eea{\end{eqnarray}}
\def\beann{\begin{eqnarray*}}
\def\eeann{\end{eqnarray*}}

\def\ns{\hspace{-1mm}}

\def\proof{\noindent{\bf{\em Proof:}\ \ }}
\def\QED{\mbox{\rule[0pt]{1.5ex}{1.5ex}}}
\def\endproof{\hspace*{\fill}~\QED\par\endtrivlist\unskip}

\newcommand{\ima}{\operatorname{im}}

\newcommand{\defi}{\stackrel{\text{\tiny def}}{=}}

\def\gD{{\cal D}}

\def\gI{{\cal I}}
\def\gJ{{\cal J}}
\def\gK{{\cal K}}

\def\gM{{\cal M}}
\def\gN{{\cal N}}

\def\gP{{\cal P}}
\def\gQ{{\cal Q}}
\def\gR{{\cal R}}
\def\gS{{\cal S}}
\def\gT{{\cal T}}
\def\gU{{\cal U}}
\def\gV{{\cal V}}
\def\gW{{\cal W}}
\def\gX{{\cal X}}
\def\gY{{\cal Y}}
\def\gZ{{\cal Z}}

\def\bmat{\left[ \begin{array}}
\def\emat{\end{array} \right]}

\def\bmat{\left[ \begin{array}}
\def\emat{\end{array} \right]}
\def\bsmat{\left[ \begin{smallmatrix}}
\def\esmat{\end{smallmatrix} \right]}

\def\l{{\lambda}}
\def\gP{{\cal P}}

\def\gU{{\cal U}}

\def\gR{{\cal R}}

\def\gV{{\cal V}}

\def\gS{{\cal S}}
\def\gK{{\cal K}}
\def\gT{{\cal T}}

\def\gX{{\cal X}}

\def\w{{w}}
\def\x{{{x}}}
\def\u{{u}}
\def\y{{y}}
\def\z{{z}}

\def\m{{m}}
\def\n{{n}}

\def\s{{s}}

\def\p{{p}}

\def\sy{\scriptscriptstyle}

\newcommand{\spanR}{\operatorname{span}}

\begin{document}
\title{\LARGE{Fixed poles of the disturbance decoupling problem by dynamic output feedback for biproper systems}}

\author{Fabrizio Padula and Lorenzo Ntogramatzidis
 \thanks{Fabrizio Padula and L. Ntogramatzidis are with the Department of Mathematics and
Statistics, Curtin University, Perth,
Australia. E-mail: {\tt \{Fabrizio.Padula,L.Ntogramatzidis\}@curtin.edu.au. }}\\                               
}

\maketitle

\vspace{-1cm}

\IEEEpeerreviewmaketitle

\begin{abstract}                         
This paper investigates the disturbance decoupling problem by dynamic output feedback in the general case of systems with possible input-output feedthrough matrices. 
In particular, we aim to extend the geometric condition based on self-boundedness and self-hiddenness, which as is well-known enables to solve the decoupling problem without requiring eigenspace computations. We show that, exactly as in the case of zero feedthrough matrices, this solution maximizes the number of assignable eigenvalues of the closed-loop. 
Since in this framework we are allowing every feedthrough matrix to be non-zero, an issue of well-posedness of the feedback interconnection arises, which affects the way the solvability conditions are structured. We show, however, that the further solvability condition which originates from the problem of well-posedness is well-behaved in the case where we express such condition in terms of self bounded and self hidden subspaces.

\end{abstract}

\section{Introduction}
\label{secintro}

One of the most articulated, rich and interesting problems within the family of disturbance decoupling problems \cite{Basile-M-92,Wonham-M-70} is, without any doubt, the disturbance decoupling problem by dynamic output feedback. 
 The first paper which provided a solution to this problem is \cite{Schumacher-80-TAC-N6}, via the introduction of pairs of subspaces known as $(C,A,B)$-pairs. Around the same time, the same problem with the additional requirement of internal stability was addressed in \cite{Willems-C-81} and \cite{Imai-A-81}.
In \cite{Basile-MP-87-JOTAI}, an alternative geometric solution for the decoupling problem with internal stability was presented, in which the solvability conditions are expressed in terms of certain self bounded and self hidden subspaces \cite{Basile-M-92}. This solution, differently from the previous ones based on stabilizability and detectability subspaces, avoids the use of eigenspaces, and only relies on subspaces obtainable from finite sequences involving only additions and intersections of subspaces, and images and counterimages of linear mappings, and therefore remains at the fundamental level of finite arithmetics.

In \cite{DMCuellar-M-01} it was shown that -- similarly to what happens for the disturbance decoupling by state feedback \cite{Malabre-MG-DMC-97} -- in the solution of the disturbance decoupling problem by dynamic output feedback there is a number of closed-loop eigenvalues that are fixed for any feedback controller that solves the decoupling problem; these unassignable eigenvalues are called {\em fixed poles} of the disturbance decoupling problem. 
Importantly, in \cite{DMCuellar-M-01} it was proved that, under very mild assumptions, the solution of the disturbance decoupling problem with dynamic output feedback based on the idea of self boundedness and self hiddenness is the best in terms of assignability of the closed-loop dynamics, see also \cite{DMCuellar-97}. The set of assumptions of \cite{DMCuellar-M-01} were then weakened further in \cite{DMCuellar-M-ECC}.
Most of the literature in geometric control has been developed for those systems which have zero feedthrough between the input and the output. 

The disturbance decoupling problem with dynamic output feedback and nonzero feedthrough has been completely solved in terms of stabilizability and detectability subspaces in \cite{Stoorvogel-W-91}. More recently, the approach based on self boundedness and self hiddenness has been generalized in \cite{N-TAC-08} for the disturbance decoupling problem with static state-feedback. In \cite{N-TAC-08}, the result of \cite{Malabre-MG-DMC-97} on the fixed poles was also generalized to biproper systems.

The disturbance decoupling problem by dynamic output feedback for biproper systems using the concepts of self boundedness
and self hiddenness is significantly more challenging, and it has been addressed only very recently \cite{PN-submitted}. The significant increase in mathematical complexity is due to the fact that an issue of well-posedness of the feedback interconnection arises when the feedthrough matrix between the control input
and the measurement output is allowed to be non-zero. It was observed in \cite{Stoorvogel-W-91} that the solvability conditions, when dealing with the problem in its full generality, need to take into account the well-posedness problem: this results in a condition that cannot be expressed as the typical subspace inclusion of most control/estimation problems for which a geometric solution is available, and therefore
$(C,A,B)$-pairs $(\gS,\gV)$, with $\gS$ and $\gV$ input containing and output nulling subspaces, respectively,
are not suitable to describe the solvability in this general case. In fact, the well-posedness issue requires the existence of a suitable gain matrix, herein denoted by $K$, that renders the feedback interconnection feasible, and therefore the concept of $(C,A,B)$-pair $(\gS,\gV)$ is substituted with the concept of solution triple $(\gS,\gV;K)$.
In \cite{PN-submitted}, the role that the well-posedness condition plays in the disturbance decoupling problem by
dynamic output feedback was investigated.  A surprising result of \cite{PN-submitted} is the fact that the well-posedness condition does not limit the solvability of the problem with the additional requirement of closed-loop stability: if the well-posedness condition is satisfied for the supremal/infimal subspaces, it is also satisfied for a special stabilizing pair of subspaces, one of which is self bounded and the other is self hidden.

This paper complements and completes the theory of the disturbance decoupling by dynamic output feedback for biproper systems developed in the recently published paper \cite{PN-submitted}. 
In particular, in Section~\ref{sec:fixed} we solve for the first time in the literature the fixed poles problem in the general case where all the feedthrough matrices are allowed to be non-zero. We establish a new fundamental result that shows that the well-posedness property is invariant with respect to {\em all} the possible pairs of self bounded and self hidden subspaces involved in the solution of the problem. 
To this end, we provide an exhaustive characterization of all the gain matrices $K$ that can be used to form solution triples.
This is the cornerstone upon which we build a theory of fixed poles for the biproper case, which addresses several delicate issues which arise from the well posedness of the feedback interconnection.
Specifically, we show that the best possible solution, in terms of assignability of the closed-loop eigenstructure, is based on self bounded and self hidden subspaces. In particular, we first prove that if a solution to the disturbance decoupling problem by dynamic output feedback exists, we can always build an alternative controller, by using self hidden and self bounded subspaces, that achieves more freedom in the assignability of the closed loop spectrum. Then, we characterize the set of self bounded and self hidden subspaces which minimizes the number of closed loop unassignable poles. Finally, we exploit the invariance of the well-posedness to show that each element of the aforementioned set can be used to form a solution triple, i.e., it provides a well-posed solution to the disturbance decoupling problem that maximizes the freedom in the assignment of the closed-loop eigenvalues.
Moreover, we provide a characterization of the minimal set of unassignable closed-loop eigenvalues, i.e., the fixed poles, when this minimal set exists, and we offer a lower and upper bound when the set of unassignable closed-loop eigenvalues does not have a minimum.
Interestingly, such upper bound is uniquely expressed as the fixed poles of a special pair of self bounded and self hidden subspaces. 

\par\noindent
{\bf Notation.} Given a vector space $\gX$, we denote by $0_{\scriptscriptstyle \scriptscriptstyle \gX}$ the origin of $\gX$. The image the kernel and the Moore-Penrose pseudoinverse of matrix $A$ are denoted by $\ima\,A$, $\ker\,A$ and $A^\dagger$, respectively. Given a linear map $A: \gX \longrightarrow \gY$ and a subspace $\gS$ of $\gY$, we define $A^{-1}\,\gS$ is for the inverse image of $\gS$ with respect to the linear map $A$, i.e., 
$A^{-1}\,\gS=\{x \in \gX\,|\,A\,x \in \gS\}$. 
When $A$ is square, we denote by $\sigma(A)$ the spectrum of $A$. If $A: \gX \longrightarrow \gY$ is a linear map and if $\gJ \subseteq \gX$, the restriction of the map $A$ to $\gJ$ is denoted by $A\,|\gJ$. If $\gX=\gY$ and $\gJ$ is $A$-invariant, the eigenstructure of $A$ restricted to $\gJ$ is denoted by $\sigma\,(A\,| \gJ)$. If $\gJ_1$ and $\gJ_2$ are $A$-invariant subspaces and $\gJ_1\,{\subseteq}\,\gJ_2$, the mapping induced by $A$ on the quotient space $\frac{\gJ_2}{\gJ_1}$ is denoted by $A\,| \frac{\gJ_2}{\gJ_1}$, and its spectrum by $\sigma\,\left(A\,\Big| \,\frac{\gJ_2}{\gJ_1}\right)$.
The symbol $\oplus$ stands for the direct sum of subspaces.  The symbol $\uplus$ denotes union with any common elements repeated.
Given a map $A: \gX \longrightarrow \gX$ and a subspace $\gS$ of $\gX$, $\langle A\,|\, \gS \rangle$ denotes the smallest $A$-invariant subspace of $\gX$ containing $\gS$ and $\langle \gS\,|\, A \rangle$ is the largest $A$-invariant subspace contained in $\gS$. 

\section{Problem Statements}
In what follows, $\mathbb{T}$ denotes $\real^+$ in the continuous time and $\mathbb{N}$ in the discrete time. 
The operator
$\gD$ denotes either the time derivative in the continuous time,
i.e., $\gD\, x(t)=\dot{x}(t)$, or the unit time shift
 in the discrete time, i.e., $\gD\, x(t)=x(t+1)$. 
We consider the system $\Sigma$ ruled by
\beann
\Sigma: \left\{\begin{array}{rcl}
\gD\,{\x}(t) \ns&\ns = \ns&\ns A\,\x(t)+B\,\u(t)+H\,\w(t) \\[0mm]
\y(t) \ns&\ns = \ns&\ns C\,\x(t)+D_y\,\u(t)+G_y\,\w(t)\\[0mm]
\z(t) \ns&\ns = \ns&\ns E\,\x(t)+D_z\,\u(t)+G_z\,\w(t),
\end{array}\right.
\eeann
where, for $t\in \mathbb{T}$, $\x(t)\in \gX=\real^n$ is the state,  $\u(t)\in \gU=\real^m$ is the control input, $\w(t)\in \gW=\real^q$ is the disturbance input, $\y(t)\in \gY=\real^p$ is the measurement output and  $\z(t)\in \gZ=\real^r$ is the to-be-controlled output. Consider the regulator
\beann
\Sigma_{\scriptscriptstyle C}: \left\{\begin{array}{rcl}
\gD{\p}(t) \ns&\ns = \ns&\ns A_c\,\p(t)+B_c\,\y(t) \\[0mm]
\u(t) \ns&\ns = \ns&\ns C_c\,\p(t)+D_c\,\y(t),
\end{array}\right.
\eeann
where, for all $t\in \mathbb{T}$, $\p(t)\in \gP=\real^s$. 
We want to control the system $\Sigma$ with the regulator $\Sigma_{\scriptscriptstyle C}$ in such a way that in the closed-loop system the output $\z$
does not depend on the disturbance input $\w$. 
We say that the feedback interconnection of system $\Sigma$ with the regulator $\Sigma_{\scriptscriptstyle C}$ is well posed if $I-D_y\,D_c$ is non-singular, see \cite[Chpt.~3]{Trentelman-SH-01}. In such case, the closed-loop system can be written in state-space form as
\bea
\label{cl}
\Sigma_{\scriptscriptstyle CL}: \left\{\begin{array}{rcl}\!\!
\gD{\hat{x}}(t) 
\ns&\ns \!=\! \ns&\ns \widehat{A}\,\hat{x}(t) +
\widehat{H}\, \w(t) \\[0mm]
\z(t)
\ns&\ns \!=\! \ns&\ns \widehat{C}\,\hat{x}(t)+ \widehat{G}\,\w(t),
\end{array}\right.
\eea
where $\hat{x}(t)=\bsmat x(t) \\[1mm] p(t) \esmat$ is the extended state, and 
\bea
\widehat{A} \ns&\ns \defi \ns&\ns \bmat{cc}
A\!+\!B D_c W C & B C_c\!+\!B D_c W D_y C_c \!\\
\!B_c W C & A_c\!+\!B_c W D_y C_c\!\!  \emat\!\!, \!\!\quad \!\!\!\widehat{H} \defi\bmat{c} \!H\!+\!B D_c W G_y \! \nonumber  \\
 \!\! B_c W G_y \!\! \emat\!\!, \nonumber \\
\widehat{C} \ns&\ns \defi \ns&\ns[\begin{array}{cc}
\!E\!+\!D_z\,D_c\,
W\,
C & D_z\,C_c\!+\!D_z\,D_c\,W\,D_y\,C_c \!\end{array}], \nonumber\\
 \widehat{G} \ns&\ns \defi  \ns&\ns G_z\!+\!D_z\,D_c\,W\,G_y,
\eea
where $W\defi (I-D_y\,D_c)^{-1}$.  Using the matrix inversion lemma\footnote{If $P$, $R$ and $P+Q\,R\,S$ are invertible matrices, we have $(P+Q\,R\,S)^{-1}=P^{-1}-P^{-1}\,Q\,(R^{-1}+S\,P^{-1}\,Q)^{-1}\,S\,P^{-1}$.
}
we can re-write the submatrix of $\widehat{A}$ in position $(1,2)$ as $B\,(I-D_c\,D_y)^{-1}\,C_c$.
We require 
that $I-D_y\,D_c$ be non-singular, i.e., that the interconnection be well-posed. 
The transfer function of the closed-loop system $\Sigma_{\scriptscriptstyle CL}$ is
$G_{\z,\w}(\l)=\widehat{C}\,(\l\,I-\widehat{A})^{-1}\widehat{H}+\widehat{G}$,
where $\l$ represents the $s$ variable of the Laplace transform in the continuous time or the $z$ variable of the $\gZ$-transform in the discrete time.

A fundamental problem in control theory is the disturbance decoupling by dynamic output feedback (DDPDOF in short), which involves finding a controller $\Sigma_{\scriptscriptstyle C}$ such that
the feedback interconnection of $\Sigma$ with $\Sigma_{\scriptscriptstyle C}$ is well posed and the transfer function matrix $G_{\z,\w}(\l)$
 of the closed-loop system $\Sigma_{\scriptscriptstyle CL}$ is zero. 
When the feedthrough matrices are different from zero, conditions for the solvability of the DDPDOF are available in the literature, see e.g. \cite{Stoorvogel-W-91}. However, these conditions, in general, overconstrain the closed-loop spectrum. In this paper, we are concerned with finding a strategy to tackle the DDPDOF that guarantees maximal freedom in the assignability of the closed-loop spectrum. We will show in Section~\ref{sec:fixed} that this problem is well posed, and it can be reformulated in geometric terms along the same lines of \cite{DMCuellar-M-01} for systems without feedthrough matrices.

\section{Geometric background}
Given $A,B,C,D$, we briefly define the main geometric concepts used in this paper. We refer the reader to \cite{Basile-M-92,Trentelman-SH-01,PN-submitted} for more details.
We denote by $\gR$ the reachable subspace of the pair $(A,B)$,  i.e., $\gR=\langle A\,|\,\ima B\rangle$, and by $\gQ$ the unobservable subspace of the pair $(C,A)$, i.e., $\gQ=\langle \ker C\,|\,A\rangle$.

A subspace $\gV$ is said to be an $(A,B,C,D)$-{\em output nulling subspace} if $\bsmat A \\[1mm]
C\esmat \,\gV\subseteq (\gV\oplus 0_{\scriptscriptstyle \scriptscriptstyle \gY})+\ima \bsmat B\\[1mm] D \esmat$, or, equivalently, if 
there exists $F$ such that 
$\bsmat A+B\,F \\[1mm] C+D\,F \esmat \,\gV\subseteq \gV\oplus 0_{\scriptscriptstyle \scriptscriptstyle \gY}$: in this case
we say that $F$ is an $(A,B,C,D)$-{\em output nulling friend} of $\gV$. We denote by $\mathfrak{F}_{\scriptscriptstyle (A,B,C,D)}(\gV)$ the set of $(A,B,C,D)$-output nulling friends of $\gV$.
 We denote by $\gV^\star_{\scriptscriptstyle (A,B,C,D)}$ the largest $(A,B,C,D)$-output nulling subspace.
Given an $(A,B,C,D)$-output nulling subspace $\gV$, the $(A,B,C,D)$-reachability subspace $\gR_{\gV}$ on $\gV$ is defined as $\gR_{\gV}=\langle A+B\,F\,|\,\gV\cap B\,\ker D\rangle$, where $\mathfrak{F}_{\scriptscriptstyle (A,B,C,D)}(\gV)$.
The {\em fixed poles} of $\gV$ are the eigenvalues of $A+B\,F$ that do not depend on $F\in \mathfrak{F}_{\scriptscriptstyle (A,B,C,D)}(\gV)$, and are given by
\bea
\label{sfixV_II0}
\sigma_{\rm fixed}(\gV) \defi \sigma \Big(A+B\,F\,\Big|\,\frac{\gX}{\gV+\gR}\Big) \uplus 
\sigma \Big(A+B\,F\,\Big|\,\frac{\gV}{\gR_{\gV}}\Big).
\eea
The remaining eigenvalues of $A+B\,F$ are freely assignable with a suitable $F\in \mathfrak{F}_{\scriptscriptstyle (A,B,C,D)}(\gV)$.
It is easy to see that $\sigma_{\rm fixed}(\gV)$ can be alternatively expressed as
\bea
\label{sfixV}
\sigma_{\rm fixed}(\gV) = \sigma \Big(A+B\,F\,\Big|\,\frac{\gX}{\gR}\Big) \uplus 
\sigma \Big(A+B\,F\,\Big|\,\frac{\gV\cap \gR}{\gR_{\gV}}\Big), 
\eea
for $F\in \mathfrak{F}_{\scriptscriptstyle (A,B,C,D)}(\gV)$, see \cite{DMCuellar-M-01,Malabre-MG-DMC-97}.
An $(A,B,C,D)$-output nulling subspace $\gR$ for which an output nulling friend $F$ exists such that the spectrum of $A+B\,F\,|\,\gR$ is arbitrary is called an $(A,B,C,D)$-{\em reachability output nulling subspace}. We denote by $\gR^\star_{\scriptscriptstyle (A,B,C,D)}$ the maximum of the set of $(A,B,C,D)$-reachability output nulling subspaces, and it coincides with the
 output nulling reachability subspace on $\gV^\star_{\scriptscriptstyle (A,B,C,D)}$, i.e., $\gR^\star_{\scriptscriptstyle (A,B,C,D)} = \gR_{\gV^\star_{\scriptscriptstyle (A,B,C,D)}}$. 
  The spectrum of $A+B\,F\,|\,\frac{\gV^\star_{\scriptscriptstyle (A,B,C,D)}}{\gR^\star_{\scriptscriptstyle (A,B,C,D)}}$ is the {\em invariant zero structure} of the system, and it is denoted by $Z_{\scriptscriptstyle (A,B,C,D)}$.
  
We say that an $(A,B,C,D)$-output nulling subspace $\gV$ is $(A,B,C,D)$-{\em self bounded} if $\gV\supseteq \gV^\star_{\scriptscriptstyle (A,B,C,D)}\cap B\,\ker D$, or, equivalently, if $\gV \supseteq \gR^\star_{\scriptscriptstyle (A,B,C,D)}$.
It follows immediately that $\gR^\star_{\scriptscriptstyle (A,B,C,D)}$ and $\gV^\star_{\scriptscriptstyle (A,B,C,D)}$ are $(A,B,C,D)$-self bounded subspaces.
Every $(A,B,C,D)$-output nulling friend of $\gV^\star_{\scriptscriptstyle (A,B,C,D)}$ is also an $(A,B,C,D)$-output nulling friend of $\gR^\star_{\scriptscriptstyle (A,B,C,D)}$.
The intersection of $(A,B,C,D)$-self bounded subspaces is $(A,B,C,D)$-self bounded. Thus, the set $\Phi_{\scriptscriptstyle (A,B,C,D)}$ of $(A,B,C,D)$-self bounded subspaces has a maximum, which is $\gV^\star_{\scriptscriptstyle (A,B,C,D)}$, and a minimum, which is $\gR^\star_{\scriptscriptstyle (A,B,C,D)}$.

A subspace $\gS$ is an $(A,B,C,D)$-{\em input containing subspace} if $[\begin{array}{cc} A & B \end{array}]\left((\gS \oplus \gU) \cap \ker [\begin{array}{cc} C & D \end{array}]\right)\subseteq \gS$, or, equivalently if there exists $G$ such that 
$[\begin{array}{cc} A+G\,C & B+G\,D \end{array}] \,(\gS\oplus \gU) \subseteq \gS$, and we say that $G$ is an $(A,B,C,D)$-{\em input containing friend} of $\gS$.
  We denote by $\mathfrak{G}_{\scriptscriptstyle (A,B,C,D)}(\gS)$ the set of $(A,B,C,D)$-input containing friends of $\gS$. 
The smallest $(A,B,C,D)$-input containing subspace is denoted by $\gS^\star_{\scriptscriptstyle (A,B,C,D)}$. There holds also $\gS^\star_{\scriptscriptstyle (A,B,C,D)}=\bigl(\gV^\star_{\scriptscriptstyle (A^\top,C^\top,B^\top,D^\top)}\bigr)^\perp$, where $(A^\top,C^\top,B^\top,D^\top)$ is the dual of $(A,B,C,D)$.
Given an $(A,B,C,D)$-input containing subspace $\gS$ and  $G\in \mathfrak{G}_{\scriptscriptstyle (A,B,C,D)}(\gS)$, we define the $(A,B,C,D)$-detectability subspace associated to it as 
$\gQ_{\gS}=\langle \gS+C^{-1}\ima D \,|\,A+G\,C
\rangle$, which is the orthogonal complement of the reachability subspace on $\gS^\perp$.
 The fixed poles of $\gS$ are the  eigenvalues of $A+G\,C$ that do not depend on $G\in \mathfrak{G}_{\scriptscriptstyle (A,B,C,D)}(\gS)$, and they can be written as
\bea
\sigma_{\rm fixed}(\gS) \defi \sigma \Big(A+G\,C\,\Big|\,\frac{\gQ_{\gS}}{\gS}\Big) \uplus 
\sigma \Big(A+G\,C\,\Big|\,\gS\cap \gQ \Big),
\eea
or, which is the same, as
$\sigma_{\rm fixed}(\gS) = \sigma \left(A+G\,C\,\Big|\,\frac{\gQ_{\gS}}{\gS+\gQ}\right) \uplus 
\sigma \left(A+G\,C\,\Big|\,\gQ\right)$ for $G\in \mathfrak{G}_{\scriptscriptstyle (A,B,C,D)}(\gS)$.
An input containing subspace $\gQ$ for which an $(A,B,C,D)$-input containing friend $G$ exists such that the spectrum of $A+G\,C\,|\,\frac{\gX}{\gQ}$ is arbitrary is called an $(A,B,C,D)$-{\em unobservability input containing subspace}. The set of $(A,B,C,D)$-unobservability input containing subspaces has a minimum, that we denote by $\gQ^\star_{\scriptscriptstyle (A,B,C,D)}$. 
There holds also $\gQ^\star_{\scriptscriptstyle (A,B,C,D)} = \gQ_{\gS^\star_{\scriptscriptstyle (A,B,C,D)}}$.  The spectrum $A+G\,C\,|\,\frac{\gQ^\star_{\scriptscriptstyle (A,B,C,D)}}{\gS^\star_{\scriptscriptstyle (A,B,C,D)}}$ coincides with the {\em invariant zero structure} of the system:
$Z_{\scriptscriptstyle (A,B,C,D)}=\sigma\bigl(A+B\,F\,|\,\frac{\gV^\star_{\scriptscriptstyle (A,B,C,D)}}{\gR^\star_{\scriptscriptstyle (A,B,C,D)}}\bigr)=\sigma\bigl(A+G\,C\,|\,\frac{\gQ^\star_{\scriptscriptstyle (A,B,C,D)}}{\gS^\star_{\scriptscriptstyle (A,B,C,D)}}\bigr)$. Finally, we recall that $\gQ^\star_{\scriptscriptstyle (A,B,C,D)}$ is the dual of $\gR^\star_{\scriptscriptstyle (A,B,C,D)}$, i.e., $\gQ^\star_{\scriptscriptstyle (A,B,C,D)}=\bigl(\gR^\star_{\scriptscriptstyle (A^\top,C^\top,B^\top,D^\top)}\bigr)^\perp$.

An $(A,B,C,D)$-input containing subspace $\gS$ is $(A,B,C,D)$-{\em self hidden} if  $\gS \subseteq \gS^\star_{\scriptscriptstyle (A,B,C,D)}+C^{-1} \ima  D$ or, equivalently, if 
 $\gS \subseteq \gQ^\star_{\scriptscriptstyle (A,B,C,D)}$.
Thus, $\gQ^\star_{\scriptscriptstyle (A,B,C,D)}$ and $\gS^\star_{\scriptscriptstyle (A,B,C,D)}$ are $(A,B,C,D)$-self hidden subspaces.
Every $(A,B,C,D)$-input containing friend of $\gS^\star_{\scriptscriptstyle (A,B,C,D)}$ is also an $(A,B,C,D)$-input containing friend of $\gQ^\star_{\scriptscriptstyle (A,B,C,D)}$, and 
the sum of $(A,B,C,D)$-self hidden subspaces is $(A,B,C,D)$-self hidden. We define $\Psi_{\scriptscriptstyle (A,B,C,D)}$ to be the set of $(A,B,C,D)$-self hidden subspaces: $\Psi_{\scriptscriptstyle (A,B,C,D)}$ has a maximum $\gQ^\star_{\scriptscriptstyle (A,B,C,D)}$ and a minimum, which is $\gS^\star_{\scriptscriptstyle (A,B,C,D)}$.

\section{Solution triples}
We begin by first presenting the following result, see \cite[Lemma 3.2]{Stoorvogel-W-91}. The proof can be carried out along the same lines of the proof of \cite[Lemma 6.3]{Trentelman-SH-01}; nevertheless, a brief proof is given here because it forms the basis of the parameterization established in Theorem~\ref{theparam}.

\begin{lemma}
\label{lemstoor}
Let $n_1,n_2,m,p\in \mathbb{N} \setminus \{0\}$. Consider a subspace
$\gM$ of $\real^{n_2}$ and a subspace $\gN$ of $\real^{n_1}$.
Consider  $\tilde{A}\in \real^{n_1 \times n_2}$, $\tilde{B}\in \real^{n_1 \times m}$ and $\tilde{C}\in \real^{p\times n_2}$.
Then, $\tilde{A}\,\gM \subseteq \gN+\ima \tilde{B}$ and $\tilde{A}\,(\gM\cap \ker \tilde{C}) \subseteq \gN$
if and only if there exists  $K\in \real^{m \times p}$ such that
$(\tilde{A}+\tilde{B}\,K\,\tilde{C})\,\gM\subseteq \gN$.
\end{lemma}

\proof
{\bf (Only if)}. Let $k=\dim(\gM \cap \ker \tilde{C})$, $r=\dim \gM$, and let
 $\m_1,\m_2,\ldots,\m_k$ be a basis of $\gM \cap \ker \tilde{C}$. Let $\m_1,\ldots,\m_r$ be a basis of $\gM$. Since $\tilde{A}\,\gM \subseteq \gN+\ima \tilde{B}$, for every $i\in \{1,\ldots,r\}$ there exist $\u_i \in \real^m$ and $\n_i\in \gN$ such that
\bea
\tilde{A}\,\m_i=\n_i+\tilde{B}\,\u_i. \label{AsvB1}
\eea
We now prove that the vectors $\tilde{C}\,\m_{k+1},\ldots,\tilde{C}\,\m_r$ are linearly independent. Consider the zero linear combination
$\alpha_{k+1}\,\tilde{C}\,\m_{k+1}+\ldots+\alpha_{r}\,\tilde{C}\,\m_{r}=0$.
The vector $\m\defi \alpha_{k+1}\,\m_{k+1}+\ldots+\alpha_{r}\,\m_{r}$ is such that $\tilde{C}\,\m=0$, so that $\m\in \gM\cap \ker \tilde{C}$. It follows that $\m$ is also a linear combination of $\m_1,\ldots,\m_k$, i.e., we can write
$\m=\alpha_1\,\m_1+\ldots+\alpha_k\,\m_k$.
Therefore
$\m-\m=\alpha_1\,\m_1+\ldots+\alpha_k\,\m_k-\alpha_{k+1}\,\m_{k+1}-\ldots-\alpha_{r}\,\m_{r}=0$,
so that $\alpha_i=0$ for all $i\in \{1,\ldots,r\}$ since $\m_1,\ldots,\m_r$ are linearly independent. In particular, $\alpha_i=0$ for all $i\in \{k+1,\ldots,r\}$.
Thus, $\tilde{C}\,\m_{k+1},\ldots,\tilde{C}\,\m_r$ are linearly independent, so that there exists $K\in \real^{m \times p}$ such that $K\,\tilde{C}\,\m_i=-\u_i$ for all $i\in \{k+1,\ldots,r\}$ since $[\begin{array}{cccc} \tilde{C}\,\m_{k+1} & \ldots & \tilde{C}\,\m_r\end{array}]$ is full column-rank. Using these identities into (\ref{AsvB1}) gives 
$(\tilde{A}+\tilde{B}\,K\,\tilde{C})\,\m_i=\n_i\in \gN$
for all $i\in \{k+1,\ldots,r\}$ and
$(\tilde{A}+\tilde{B}\,K\,\tilde{C})\,\m_i=\tilde{A}\,\m_i \in \tilde{A}\,(\gM\cap \ker \tilde{C})\subseteq\gN$
for all $i\in \{1,\ldots,k\}$. Thus, $(\tilde{A}+\tilde{B}\,K\,\tilde{C})\,\m_i\in \gN$ for all $i\in \{1,\ldots,r\}$. Since $\m_1,\ldots,\m_r$ span $\gM$, then $(\tilde{A}+\tilde{B}\,K\,\tilde{C})\,\gM\subseteq \gN$.

\noindent {\bf (If)}. It is easy to see that $(\tilde{A}+\tilde{B}\,K\,\tilde{C})\,\gM\subseteq \gN$ implies $\tilde{A}\,\gM \subseteq \gN+\ima \tilde{B}$ and $\tilde{A}\,(\gM\cap \ker \tilde{C}) \subseteq \gN$.
\endproof


The following result uses the proof of the latter to determine an exhaustive parameterization of the set of all matrices $K\in \real^{m \times p}$ such that $(\tilde{A}+\tilde{B}\,K\,\tilde{C})\gM\subseteq \gN$.


\begin{theorem}
\label{theparam}
Let $N$ be a basis matrix of $\gN$. Let $\m_1,\ldots,\m_r$ be a basis of $\gM$ such that  $\m_1,\m_2,\ldots,\m_k$ is a basis of $\gM\cap \ker \tilde{C}$, and let $M= [\begin{array}{ccccc} \m_{k+1} &\ldots&\m_{r} \end{array}]$.
The set of matrices $K\in \real^{m \times p}$ such that
$(\tilde{A}+\tilde{B}\,K\,\tilde{C})\,\gM\subseteq \gN$ is parameterized in terms of the  matrices $H_1$ and $H_2$ (of suitable sizes) as
\bea
\label{param}
K=-R_2\,\tilde{A}\,M\,(\tilde{C}\,M)^\dagger-\Phi_2\,H_1\,(\tilde{C}\,M)^\dagger+H_2\,\Psi,
\eea
where
\begin{itemize}
\item $\bsmat R_1 \\[1mm] R_2\esmat=[\begin{array}{cc} N & \tilde{B} \end{array}]^\dagger$, partitioned conformably;
\item $\Phi=\bsmat \Phi_1\\[1mm] \Phi_2\esmat$ is a basis matrix of $\ker [\begin{array}{cc} N & \tilde{B} \end{array}]$, partitioned conformably;
\item $\Psi$ is a basis matrix of the left null-space of $\tilde{C}\,M$, i.e., $\Psi$ is full row-rank and $\Psi\,\tilde{C}\,M=0$.
\end{itemize}
\end{theorem}
\proof
Consider the {\em only if} part of the proof of Lemma~\ref{lemstoor}.
We rewrite (\ref{AsvB1}) as $\tilde{A}\,\m_i=N\,\xi_i+\tilde{B}\,\u_i$, whose set of solutions is parameterized in $h_i$ as
$\bsmat \xi_i \\[1mm] \u_i \esmat=[\begin{array}{cc} N & \tilde{B} \end{array}]^\dagger\,\tilde{A}\,\m_i+\Phi\,h_i$ for all $i\in\{k+1,\ldots,r\}$. With our definitions, with $\Xi=[\begin{array}{ccc} \xi_{k+1} & \ldots & \xi_r \end{array}]$ and 
$U=[\begin{array}{ccc} \u_{k+1} & \ldots & \u_r \end{array}]$, such equation can be re-written as $\bsmat \Xi \\[1mm] U \esmat=\bsmat R_1 \\[1mm] R_2\esmat\,\tilde{A}\,M+\Phi\,H_1$. Now,
$K\in \real^{m \times p}$ must satisfy $K\,\tilde{C}\,\m_i=-\u_i$ for all $i\in \{k+1,\ldots,r\}$, i.e., $K\,\tilde{C}\,M=-U=-R_2\,\tilde{A}\,M-\Phi\,H_1$, from which (\ref{param}) immediately follows. 
We now show that $K$ does not depend on the chosen basis of $\gM$. Indeed, let $\x\in \gM$. Let $\x=c_1\,\m_1+\ldots+c_r\,\m_r$. By linearity, from $\tilde{A}\,\m_1-\tilde{B}\,\u_1\in \gN$, $\tilde{A}\,\m_2-\tilde{B}\,\u_2\in \gN$, $\ldots$, $\tilde{A}\,\m_r-\tilde{B}\,\u_r\in \gN$, we also have
$\tilde{A}\,\x-\tilde{B}\,(c_1\,\u_1+\ldots+c_r\,\u_r)\in \gN$.
We only need to prove that $K\,\tilde{C}\,\x=-c_1\,\u_1-\ldots-c_r\,\u_r$. Indeed,
$K\,\tilde{C}\,(c_1\,\m_1+\ldots+c_r\,\m_r)=-c_1\,\u_1-\ldots-c_r\,\u_r$
follows from $K\,\tilde{C}\,\m_i=-\u_i$. If we replace the set $\{\m_1,\ldots,\m_r\}$ as a basis for $\gM$ adapted to $\gM\cap \ker \tilde{C}$ with another linearly independent set $\{\m_1,\ldots,\m_{r-1},\x\}$,  matrix $K$ is the same, i.e., it still satisfies $K\,\tilde{C}\,\m_i = -\u_i$ where $\tilde{A}\,\m_i-\tilde{B}\,\u_i\in \gN$ for $i\in\{k+1,\ldots,r-1\}$ and for $i=r$ we have $K\,\tilde{C}\,\x=-c_1\,\u_1-\ldots-c_r\,\u_r$, where $\tilde{A}\,\x-\tilde{B}\,(c_1\,\u_1+\ldots+c_r\,\u_r)\in \gN$.
It remains to prove that this parameterization is exhaustive. 
Consider any basis of $\gM$ adapted to $\gM\cap \ker \tilde{C}$. We need to show that a matrix $K$ such that $(\tilde{A}+\tilde{B}\,K\,\tilde{C})\,\gM\subseteq \gN$ must satisfy
$K\,\tilde{C}\,\m_i=-\u_i$ for $i\in\{k+1,\ldots,r\}$ for some $\u_i$ such that $\tilde{A}\,\m_i-\tilde{B}\,\u_i\in \gN$. Rewriting the latter as $
\tilde{A}\,\m_i+\tilde{B}\,K\,\tilde{C}\,\m_i \in \gN$, the statement follows from $(\tilde{A}+\tilde{B}\,K\,\tilde{C})\,\gM\subseteq \gN$.
\endproof

The following two results have been proved in \cite{PN-submitted}.

\begin{lemma}
\label{lem3cond}
Let $\gV$ be $(A,B,E,D_z)$-output nulling and let $\gS$ be  $(A,H,C,G_y)$-input containing. If
\begin{description}
\item{\bf (a)} $\ima \bmat{c} H \\ G_z \emat\subseteq (\gV \oplus 0_{\scriptscriptstyle \gZ}) +\ima \bmat{c} B \\ D_z\emat$;
\item{\bf (b)} $\ker \,[\begin{array}{cc} E & G_z \end{array}] \supseteq (\gS \oplus \gW)  \cap \ker \,[\begin{array}{cc}  C & G_y \end{array}] $;
\item{\bf (c)} $\gS\subseteq \gV$
\end{description}
then there exists an output feedback matrix $K$ such that 
 \bea
 \label{clK}
 \bmat{cc} A+B\,K\,C & H+B\,K\,G_y \\
 E+D_z\,K\,C & G_z+ D_z\,K\,G_y \emat(\gS \oplus \gW)\subseteq \gV \oplus 0_{\scriptscriptstyle \gZ}.
 \eea
 Conversely, if $K$ exists such that (\ref{clK}) holds, then {\bf (a-b)} hold. 
\end{lemma}


\begin{theorem}
\label{the66ths}
DDPDOF is solvable if and only if there exist
an $(A,B,E,D_z)$-output nulling subspace $\gV$, an $(A,H,C,G_y)$-input containing subspace $\gS$ and a matrix $K\in \real^{m \times p}$ such that

\begin{description}
\item{\bf (i)} $\ima \bmat{c} H \\[-1mm] G_z \emat\subseteq (\gV \oplus {0}_{\gZ}) +\ima \bmat{c} B \\[-1mm]  D_z\emat$;
\item{\bf (ii)} $\ker \,[\begin{array}{cc} E & G_z \end{array}] \supseteq (\gS \oplus \gW)  \cap \ker \,[\begin{array}{cc}  C & G_y \end{array}] $;
\item{\bf (iii)} $\gS\subseteq \gV$;
\item{\bf (iv)} $I+K\,D_y$ is non-singular, and $K$ satisfies (\ref{clK}).
\end{description}
\end{theorem}

When DDPDOF
 is solvable, by Theorem~\ref{the66ths} there exist an $(A,B,E,D_z)$-output nulling subspace $\gV$, an 
$(A,H,C,G_y)$-input containing subspace $\gS$ and a matrix $K$ such that {\bf (i-iv)} are satisfied, and in this case we say that $(\gS,\gV;K)$ is a {\em solution triple} for DDPDOF. 
Solution triples generalize the notion of $(C,A,B)$-pairs of systems with zero feedthrough matrices.
Theorem~\ref{the66ths} does not offer a method for computing the subspaces $\gV$ and $\gS$, but, if a solution to DDPDOF exists, clearly it can always been obtained using $\gV^\star_{\scriptscriptstyle (A,B,E,D_z)}$ and $\gS^\star_{\scriptscriptstyle (A,H,C,G_y)}$ in place of $\gV$ and $\gS$, respectively, \cite{Stoorvogel-W-91}, so that, in particular, $\gV^\star_{\scriptscriptstyle (A,B,E,D_z)}$ and $\gS^\star_{\scriptscriptstyle (A,H,C,G_y)}$ satisfy 
\bea
 \ima \bmat{c} H \\[-1mm] G_z \emat \ns&\ns \subseteq  \ns&\ns (\gV^\star_{\scriptscriptstyle (A,B,E,D_z)} \oplus 0_{\scriptscriptstyle \gZ})+\ima \bmat{c} B\\[-1mm] D_z \emat, \label{condalpha} \\[0mm]
 \ker [\begin{array}{cc} E & G_z \end{array}] \ns&\ns \supseteq  \ns&\ns (\gS^\star_{\scriptscriptstyle (A,H,C,G_y)} \oplus \gW) \cap \ker[\begin{array}{cc}C & G_y \end{array}].\label{condbeta}
\eea
However, this solution is the one which constraints the maximum number of closed-loop eigenvalues. In this paper, as aforementioned, we are interested in a solution which maximizes the freedom in the assignment of the closed-loop spectrum.
Following the notation of \cite{Basile-M-92}, we define
\beann
\gV_m &\defi &\gR^\star_{\scriptscriptstyle (A,[\,B \;\; H \,],E,[\, D_z \;\; G_z \,])} =\min \Phi_{\scriptscriptstyle (A,[\,B \;\; H \,],E,[\, D_z \;\; G_z \,])}, \\
\gS_M &\defi &\gQ^\star_{\scriptscriptstyle \left(A,H,\bsmat {\scriptscriptstyle C}\\[0.3mm] {\scriptscriptstyle E} \esmat,\bsmat {\scriptscriptstyle G_y} \\[-0.5mm] {\scriptscriptstyle G_z} \esmat\right)}=\max \Psi_{\scriptscriptstyle \left(A,H,\bsmat {\scriptscriptstyle C}\\[0.3mm] {\scriptscriptstyle E} \esmat,\bsmat {\scriptscriptstyle G_y} \\[-0.5mm] {\scriptscriptstyle G_z} \esmat\right)}.
\eeann
If $\ima \bsmat H \\[1mm] G_z \esmat\subseteq (\gV^\star_{\scriptscriptstyle (A,B,E,D_z)}\oplus 0_{\scriptscriptstyle \gZ})+\ima  \bsmat B \\[1mm] D_z \esmat$, we have $\gV_m=\gV^\star_{\scriptscriptstyle (A,B,E,D_z)} \cap \gS^\star_{\scriptscriptstyle (A,[\,B \;\; H \,],E,[\, D_z \;\; G_z \,])}$, see Theorem~\ref{the14}. The following results have been proved in \cite{PN-submitted}.

\begin{lemma}
\label{lem22}
The subspace $\gV_m+\gS_M$ is $(A,[\begin{array}{cc} B & H \end{array}],E,[\begin{array}{cc} D_z & G_z\end{array}])$-self bounded. Moreover, if $\ima \bsmat H \\[1mm] G_z \esmat \subseteq (\gV^\star_{\scriptscriptstyle (A,B,E,D_z)} \oplus 0_{\scriptscriptstyle \scriptscriptstyle \gZ})+\ima \bsmat B\\[1mm] D_z \esmat$, 
then $\gV_m+\gS_M$  is also $(A,B,E,D_z)$-self bounded.
Dually, the subspace $\gV_m\cap \gS_M$ is $\left(A,H ,\bsmat C\\[1mm] E\esmat,\bsmat G_y \\[0.5mm] G_z\esmat\right)$-self hidden. Moreover, if $\ker [\begin{array}{cc} E & G_z \end{array}]\supseteq (\gS^\star_{\scriptscriptstyle (A,H,C,G_y)} \oplus \gW) \cap \ker[\begin{array}{cc}C & G_y \end{array}]$, then $\gV_m\cap \gS_M$ is also $(A,H,C,G_y)$-self hidden.
\end{lemma}

The next results further characterize $\gV_m$ and $\gS_M$ in the biproper case. The proof is straightforward consequence of Theorems~\ref{the12} and \ref{the13}. 

\begin{lemma}
\label{lemA1}
Let (\ref{condalpha}-\ref{condbeta}) hold. Then:
\begin{itemize}
\item Let $S_M$ be a basis matrix for $\gS_M=\gQ^\star_{\scriptscriptstyle \left(A,H,\bsmat \scriptscriptstyle{C} \\ \scriptscriptstyle{E} \esmat,\bsmat \scriptscriptstyle{G_y} \\[-0.5mm] \scriptscriptstyle{G_z} \esmat\right)}$. Then
\bea
\gS^\star_{\scriptscriptstyle (A,H,C,G_y)}  \subseteq \gV^\star_{\scriptscriptstyle (A,B,E,D_z)} \!\! \quad \Leftrightarrow \!\! \quad
\gS_M\subseteq \gR^\star_{\scriptscriptstyle (A,[\,B \;\; S_M \,],E,[\, D_z \;\; 0 \,])};
\eea
\item Let $T_m$ be a full row-rank matrix such that $\ker T_m=\gV_m=\gR^\star_{\scriptscriptstyle (A,[\,B\;\;H\,],E,[\,D_z\;\; G_z\,])}$. Then
\bea
\gS^\star_{\scriptscriptstyle (A,H,C,G_y)}  \subseteq \gV^\star_{\scriptscriptstyle (A,B,E,D_z)}  \quad  \!\! \Leftrightarrow  \!\! \quad
\gV_m\supseteq \gQ^\star_{\scriptscriptstyle \left(A,H,\bsmat \scriptscriptstyle{C} \\ \scriptscriptstyle{T_m} \esmat,\bsmat \scriptscriptstyle{G_y} \\[-0.0mm] \scriptscriptstyle{0} \esmat\right)}.
\eea
\end{itemize}
\end{lemma}
\proof
We begin proving the first point.
$(\Rightarrow)$: let $\gS^\star_{\scriptscriptstyle (A,H,C,G_y)} \subseteq \gV^\star_{\scriptscriptstyle (A,B,E,D_z)}$; then $\gV_m+\gS_M$ is 
$(A,B,E,D_z)$-self bounded. Thus, $\gV_m+\gS_M$ is contained in $\gV^\star_{\scriptscriptstyle (A,B,E,D_z)}$. In particular, also $\gS_M\subseteq \gV^\star_{\scriptscriptstyle (A,B,E,D_z)}$. Theorem~\ref{the13} ensures that $\gR^\star_{\scriptscriptstyle (A,[\,B \;\; S_M \,],E,[\, D_z \;\; 0 \,])}$ is the smallest of all
 $(A,B,E,D_z)$-self bounded subspaces containing $\gS_M$. 
 
\noindent $(\Leftarrow)$: let $\gS_M\subseteq\gR^\star_{\scriptscriptstyle (A,[\,B \;\; S_M \,],E,[\, D_z \;\; 0 \,])}$. From Theorem \ref{the12}, if $L_1$ is a basis matrix of $\gS_M$, we find $\gS_M \subseteq \gV^\star_{\scriptscriptstyle (A,B,E,D_z)}$. From $\gS_M \supseteq \gS^\star_{\scriptscriptstyle (A,H,C,G_y)}$ we obtain $\gS^\star_{\scriptscriptstyle (A,H,C,G_y)} \subseteq \gV^\star_{\scriptscriptstyle (A,B,E,D_z)}$.

 The second point can be proved by duality. 
 We rewrite the first point with the substitutions$A\rightarrow A^\top$, $B \rightarrow C^\top$, $H \rightarrow E^\top$, $C \rightarrow B^\top$, $G_y \rightarrow D_z^\top$, $E \rightarrow H^\top$, $D_z \rightarrow G_y^\top$, $G_z \rightarrow G_z^\top$. 
 The conditions (\ref{condalpha}-\ref{condbeta}) become
 \[
\ima \bsmat E^\top \\[1mm] G_z^\top \esmat \subseteq (\gV^\star_{\scriptscriptstyle (A^\top,C^\top,H^\top,G_y^\top)} \oplus 0_{\scriptscriptstyle \gW})+\ima \bsmat C^\top \\[1mm] G_y^\top \esmat
\]
 and 
 \[
 \ker [\begin{array}{cc} H^\top & G_z^\top \end{array}]\supseteq \left((\gS^\star_{\scriptscriptstyle (A^\top,E^\top,B^\top,D_z^\top)} \oplus \gZ) \cap \ker[\begin{array}{cc}B^\top & D_z^\top \end{array}]\right).
 \]
 Taking the orthogonal complement of both inclusions yields
$\ker [\begin{array}{cc} E & G_z\end{array}]\supseteq (\gS^\star_{\scriptscriptstyle (A,H,C,G_y)} \oplus \gW)\cap \ker 
 [\begin{array}{cc} C & G_y\end{array}]$ and $\ima \bsmat H \\[1mm] G_z \esmat \subseteq (\gV^\star_{\scriptscriptstyle (A,B,E,D_z)}\oplus 0_{\scriptscriptstyle \gZ})+\ima \bsmat B \\[1mm] D_z\esmat$.
 With these substitutions, the statement of this first point becomes
 \beann
&&\gS^\star_{\scriptscriptstyle (A^\top,E^\top,B^\top,D_z^\top)}  \subseteq \gV^\star_{\scriptscriptstyle (A^\top,C^\top,H^\top,G_y^\top)} \\
&& \qquad \Leftrightarrow \qquad
\gQ^\star_{\scriptscriptstyle \left(A^\top,E^\top,\bsmat B^\top \\ H^\top \esmat,\bsmat D_z^\top \\ G_z^\top \esmat\right)} \subseteq \gR^\star_{\scriptscriptstyle (A^\top,[\,C^\top \;\; N \,],H^\top,[\, G_y^\top \;\; 0 \,])},
\eeann
 where $N$ is a basis matrix of $\gQ^\star_{\scriptscriptstyle \left(A^\top,E^\top,\bsmat {\scriptscriptstyle B^\top} \\[1mm] {\scriptscriptstyle H^\top} \esmat,\bsmat {\scriptscriptstyle D_z^\top} \\ {\scriptscriptstyle G_z^\top} \esmat\right)}$. Taking the orthogonal complements of both predicates in the equivalence relation  proves that 
$\gV^\star_{\scriptscriptstyle (A,B,E,D_z)} \supseteq \gS^\star_{\scriptscriptstyle (A,H,C,G_y)}$
is equivalent to $(\gV_m=){\gR^\star_{\scriptscriptstyle (A,[\,B\;\;H\,],E,[\,D_z\;\; G_z\,])}}\supseteq 
 \gQ^\star_{\scriptscriptstyle \left(A,H,\bsmat C\\ T_m \esmat,\bsmat G_y \\ 0 \esmat\right)}$,
 where $T_m=N^\top$ is such that $\ker T_m=\gR^\star_{\scriptscriptstyle (A,[\,B\;\;H\,],E,[\,D_z\;\; G_z\,])}$ as required.
\endproof

\begin{corollary}
\label{cor8}
Let (\ref{condalpha}-\ref{condbeta}) hold.
The following are equivalent:
\begin{description}
\item{(i)} $\gS^\star_{\scriptscriptstyle (A,H,C,G_y)} \subseteq \gV^\star_{\scriptscriptstyle (A,B,E,D_z)}$;
\item{(ii)}  $\gS_M \subseteq \gV^\star_{\scriptscriptstyle (A,B,E,D_z)}$;
\item{(iii)}  $\gV_m\supseteq \gS^\star_{\scriptscriptstyle (A,H,C,G_y)}$.
\end{description}
 \end{corollary}
\proof {\em (i)} $\Leftrightarrow$ {\em (ii)}.
 If $\gS^\star_{\scriptscriptstyle (A,H,C,G_y)} \subseteq \gV^\star_{\scriptscriptstyle (A,B,E,D_z)} $, $\gV_m+\gS_M$ is $(A,B,E,D_z)$-self bounded by Lemma \ref{lem22}, so that $\gV_m+\gS_M\subseteq\gV^\star_{\scriptscriptstyle (A,B,E,D_z)}$. Thus, $\gS_M\subseteq \gV^\star_{\scriptscriptstyle (A,B,E,D_z)} $. 
  If $\gS_M \subseteq \gV^\star_{\scriptscriptstyle (A,B,E,D_z)}$, by Theorem~\ref{the12}  $\gS_M\subseteq \gR^\star_{\scriptscriptstyle (A,[\,B \;\; S_M \,],E,[\, D_z \;\; 0 \,]}$ and, by Lemma \ref{lemA1}, $\gS^\star_{\scriptscriptstyle (A,H,C,G_y)} \subseteq \gV^\star_{\scriptscriptstyle (A,B,E,D_z)} $. 
{\em (i)} $\Leftrightarrow$ {\em (iii)} can be proved by duality. 
\endproof

The following result extends \cite[Property A.4]{DMCuellar-M-01} to the biproper case.

 \begin{property}
 \label{pA4}
 Let conditions (\ref{condalpha}-\ref{condbeta}) hold.
Let $\gS^\star_{\scriptscriptstyle (A,H,C,G_y)} \subseteq \gV^\star_{\scriptscriptstyle (A,B,E,D_z)}$. Let $\gS$ be $(A,H,C,G_y)$-input containing and such that $\gV_m\cap \gS_M \subseteq \gS \subseteq \gS_M$. Then
$\gS+\gV_m\in \Phi_{\scriptscriptstyle (A,B,E,D_z)}$.
\end{property}

\proof From Corollary~\ref{cor8}, we have $\gV_m\supseteq \gS^\star_{\scriptscriptstyle (A,H,C,G_y)}$. Moreover, from Theorem~\ref{the11d}  using $\ker M=\gV_m$ and the quadruple $(A,H,C,G_y)$ in place of $(A,B,C,D)$, we find
\bea
\label{eee1}
\gS^\star_{\scriptscriptstyle (A,H,C,G_y)}=\gS^\star_{\scriptscriptstyle \left(A,H,\bsmat {\scriptscriptstyle C} \\ {\scriptscriptstyle M} \esmat,\bsmat {\scriptscriptstyle G_y} \\ {\scriptscriptstyle 0}\esmat\right)}\quad \subseteq \gQ^\star_{\scriptscriptstyle \left(A,H,\bsmat {\scriptscriptstyle C} \\ {\scriptscriptstyle M} \esmat,\bsmat {\scriptscriptstyle G_y} \\ {\scriptscriptstyle 0} \esmat\right)}.
\eea
From the second point of Lemma~\ref{lemA1} we also have $\gQ^\star_{\scriptscriptstyle \left(A,H,\bsmat {\scriptscriptstyle C} \\ {\scriptscriptstyle M} \esmat,\bsmat {\scriptscriptstyle G_y} \\ {\scriptscriptstyle 0} \esmat\right)}\subseteq \gV_m$ so that, in particular, $\gS^\star_{\scriptscriptstyle (A,H,C,G_y)} \subseteq \gV_m$.
On the other hand, using Theorem~\ref{the11d1}, with the same substitution of quadruples, we obtain  
\bea
\label{eee2}
\gS^\star_{\scriptscriptstyle (A,H,C,G_y)}=\gS^\star_{\scriptscriptstyle \left(A,H,\bsmat \scriptscriptstyle{C} \\ \scriptscriptstyle{E} \esmat,\bsmat \scriptscriptstyle{G_y} \\[-0.5mm] \scriptscriptstyle{G_z} \esmat\right)}\quad \subseteq \gQ^\star_{\scriptscriptstyle \left(A,H,\bsmat \scriptscriptstyle{C} \\ \scriptscriptstyle{E} \esmat,\bsmat \scriptscriptstyle{G_y} \\[-0.5mm] \scriptscriptstyle{G_z} \esmat\right)}.
\eea
From (\ref{eee1}) and (\ref{eee2}) we find
\bea
\label{eee3}
\gS^\star_{\scriptscriptstyle \left(A,H,\bsmat \scriptscriptstyle{C} \\ \scriptscriptstyle{E} \esmat,\bsmat \scriptscriptstyle{G_y} \\[-0.5mm] \scriptscriptstyle{G_z} \esmat\right)}
\subseteq \gQ^\star_{\scriptscriptstyle \left(A,H,\bsmat {\scriptscriptstyle C}\\ {\scriptscriptstyle M} \esmat,\bsmat {\scriptscriptstyle G_y} \\ {\scriptscriptstyle 0} \esmat\right)}.
\eea
We have
\bea
\label{eee4}
\gS \subseteq \gS_M=\gV^\star_{\scriptscriptstyle \left(A,H,\bsmat \scriptscriptstyle{C} \\ \scriptscriptstyle{E} \esmat,\bsmat \scriptscriptstyle{G_y} \\[-0.5mm] \scriptscriptstyle{G_z} \esmat\right)}+\gS^\star_{\scriptscriptstyle \left(A,H,\bsmat \scriptscriptstyle{C} \\ \scriptscriptstyle{E} \esmat,\bsmat \scriptscriptstyle{G_y} \\[-0.5mm] \scriptscriptstyle{G_z} \esmat\right)}, 
\eea
and 
$\gS\supseteq \gS^\star_{\scriptscriptstyle (A,H,C,G_y)}= \gS^\star_{\scriptscriptstyle \left(A,H,\bsmat \scriptscriptstyle{C} \\ \scriptscriptstyle{E} \esmat,\bsmat \scriptscriptstyle{G_y} \\[-0.5mm] \scriptscriptstyle{G_z} \esmat\right)}$.
Since $\gS_M\cap \gV_m$ is $(A,H,C,G_y)$-self hidden  from Lemma~\ref{lem22}, it is also input containing, and thus 
$\gS^\star_{\scriptscriptstyle \left(A,H,\bsmat \scriptscriptstyle{C} \\ \scriptscriptstyle{E} \esmat,\bsmat \scriptscriptstyle{G_y} \\[-0.5mm] \scriptscriptstyle{G_z} \esmat\right)}=
\gS^\star_{\scriptscriptstyle (A,H,C,G_y)}\subseteq \gS_M\cap \gV_m$.
Using the modular rule \cite{Trentelman-SH-01} we find
\bea
\gS+\gV_m & = & (\gS\cap \gS_M)+\gV_m 
\nonumber \\
 & = &  \left[\gS\cap \left(\gV^\star_{\scriptscriptstyle \left(A,H,\bsmat \scriptscriptstyle{C} \\ \scriptscriptstyle{E} \esmat,\bsmat \scriptscriptstyle{G_y} \\[-0.5mm] \scriptscriptstyle{G_z} \esmat\right)}+\gS^\star_{\scriptscriptstyle \left(A,H,\bsmat \scriptscriptstyle{C} \\ \scriptscriptstyle{E} \esmat,\bsmat \scriptscriptstyle{G_y} \\[-0.5mm] \scriptscriptstyle{G_z} \esmat\right)}\right)\right]+\gV_m
 \nonumber  \\
 & = &  \left[\!\!\left(\!\!\gS \! \cap \gV^\star_{\scriptscriptstyle \! \left(\!A,H,\bsmat \scriptscriptstyle{C} \\ \scriptscriptstyle{E} \esmat,\bsmat \scriptscriptstyle{G_y} \\[-0.5mm] \scriptscriptstyle{G_z} \esmat\!\!\right)}\!\!\right)\!\!+\!\!\left(\!\!\gS\cap\gS^\star_{\scriptscriptstyle  \!  \left( \! A,H,\bsmat \scriptscriptstyle{C} \\ \scriptscriptstyle{E} \esmat,\bsmat \scriptscriptstyle{G_y} \\[-0.5mm] \scriptscriptstyle{G_z} \esmat\right)}\!\right)\!\!\right]\!\!+\!\gV_m 
 \nonumber  \\
 & = &  \left(\gS\cap \gV^\star_{\scriptscriptstyle \left(A,H,\bsmat \scriptscriptstyle{C} \\ \scriptscriptstyle{E} \esmat,\bsmat \scriptscriptstyle{G_y} \\[-0.5mm] \scriptscriptstyle{G_z} \esmat\right)}\right)+\gS^\star_{\scriptscriptstyle \left(A,H,\bsmat \scriptscriptstyle{C} \\ \scriptscriptstyle{E} \esmat,\bsmat \scriptscriptstyle{G_y} \\[-0.5mm] \scriptscriptstyle{G_z} \esmat\right)}+\gV_m 
   \nonumber \\
 & = &  \left(\gS\cap \gV^\star_{\scriptscriptstyle \left(A,H,\bsmat \scriptscriptstyle{C} \\ \scriptscriptstyle{E} \esmat,\bsmat \scriptscriptstyle{G_y} \\[-0.5mm] \scriptscriptstyle{G_z} \esmat\right)}\right)+\gV_m, 
 \label{equltima}  
\eea
where we have used the inclusions $\gS\subseteq \gS_M$, 
$\gS\supseteq \gS^\star_{\scriptscriptstyle \left(A,H,\bsmat \scriptscriptstyle{C} \\ \scriptscriptstyle{E} \esmat,\bsmat \scriptscriptstyle{G_y} \\[-0.5mm] \scriptscriptstyle{G_z} \esmat\right)}$ and $\gS^\star_{\scriptscriptstyle \left(A,H,\bsmat \scriptscriptstyle{C} \\ \scriptscriptstyle{E} \esmat,\bsmat \scriptscriptstyle{G_y} \\[-0.5mm] \scriptscriptstyle{G_z} \esmat\right)} \subseteq \gV_m$, respectively.
We now show that $\gS+\gV_m$ is $(A,B,E,D_z)$-output nulling. We know that $\gV_m$ is $(A,B,E,D_z)$-self bounded and satisfies
$\ima \bsmat H \\[1mm] G_z \esmat \subseteq (\gV_m  \oplus 0_{\scriptscriptstyle \gZ})+\ima \bsmat B\\[1mm] D_z \esmat$,
so that 
$\bsmat A && H \\[1mm] E && G_z \esmat\,(\gV_m\oplus \gW)\subseteq (\gV_m\oplus 0_{\scriptscriptstyle \gZ})+\ima \bsmat B\\[1mm] D_z \esmat$.
We also know that 
$\bsmat A && H \\[1mm] E && G_z \esmat\,\left((\gS \oplus \gW)\cap \ker [\begin{array}{cc} C & G_y \end{array}]\right) \subseteq \gS\oplus 0_{\scriptscriptstyle \gZ}$,
since $[\begin{array}{cc} A & H \end{array}]\left((\gS \oplus \gW)\cap \ker [\begin{array}{cc} C & G_y \end{array}]\right) \subseteq \gS$, which comes from the definition of $(A,H,C,G_y)$-input containing subspace, and
$(\gS \oplus \gW)\cap \ker [\begin{array}{cc} C & G_y \end{array}] \subseteq
(\gS^\star_{\scriptscriptstyle (A,H,C,G_y)} \oplus \gW)\cap \ker [\begin{array}{cc} C & G_y \end{array}]\subseteq \ker [\begin{array}{cc} E & G_z \end{array}]$ from condition (\ref{condbeta}).
Thus
\beann
&&\bmat{cc} A & H \\[-1.3mm] E & G_z \emat\,\left((\gV_m\oplus \gW)+\bigl(\gS \oplus \gW)\cap \ker [\begin{array}{cc} C & G_y \end{array}]\bigr)\right)\\
&&\qquad  \qquad \subseteq \bigl((\gV_m+\gS)\oplus 0_{\scriptscriptstyle \gZ}\bigr)+\ima \bmat{c} B\\[-1.3mm] D_z \emat.
\eeann
If we prove that
\bea
\label{eq:vaffa}
(\gS\!+\!\gV_m)\oplus\gW\subseteq (\gV_m\!\oplus \!\gW)\!+\!\bigl(\gS \!\oplus\! \gW)\cap \ker [\begin{array}{cc} C & G_y \end{array}],
\eea
then it is obvious that
\bea
\bmat{c}  A   \\[-1.3mm]  \!\!  E \!\!  \emat \!\!  (\gS\!+\!\gV_m) \ns&\ns \!=\! \ns&\ns \bmat{cc} A & H \\[-1.3mm]  E & G_z \emat \!\!  (\gS+\gV_m)\oplus 0_{\scriptscriptstyle \gZ} \nonumber \\[0mm]
\ns&\ns\! \subseteq\!\ns&\ns \bmat{cc} A & H \\[-1.3mm]  E & G_z \emat\!\! \bigl((\gS+\gV_m)\oplus\gW\bigr)  \nonumber \\[0mm]
\ns&\ns \!\!\subseteq\! \ns&\ns \bmat{cc} A \!& H \\[-1.3mm]  E \!& G_z \emat\!\! \!\left(\!\!\gV_m\!\oplus \!\gW\!\!+\!\!\bigl(\gS \!\oplus \!\gW)\!\cap \ker [\begin{array}{cc} C & G_y\! \end{array}]\bigr)\!\!\right) \nonumber  \\[0mm]
\ns&\ns\!\subseteq\!\ns&\ns \bigl((\gV_m+\gS)\oplus 0_{\scriptscriptstyle \gZ}\bigr)\!+\ima \bmat{cc} B\\[-1.3mm]  D_z \emat 
\eea
says that $\gS+\gV_m$ is $(A,B,E,D_z)$-output nulling.
To prove (\ref{eq:vaffa}), we first notice that
$(\gV_m\oplus \gW)+\bigl(\gS \oplus \gW)\cap \ker [\begin{array}{cc} C & G_y \end{array}]\bigr)
=\gV_m\oplus \gW+(\gS\cap C^{-1}\,\ima G_y)\oplus \gW$.
Indeed, any $\bsmat \x \\[1mm] \w\esmat\in (\gV_m\oplus \gW)+\bigl((\gS \oplus \gW)\cap \ker [\begin{array}{cc} C & G_y \end{array}]\bigr)$ can be written as $\bsmat \x \\[1mm] \w\esmat=\bsmat \x_1 \\[1mm] \w_1\esmat+\bsmat \x_2 \\[1mm] \w_2\esmat$
with $\x_1\in \gV_m$, $\x_2\in \gS$ and $C\,\x_2+G_y\,\w_2=0$. This implies $\x_2\in C^{-1}\,\ima G_y$. Thus,
$\bsmat \x \\[1mm] \w\esmat\in \gV_m\oplus \gW+(\gS\cap C^{-1}\,\ima G_y)\oplus \gW$.
Conversely, let $\bsmat \x \\[1mm] \w\esmat\in \gV_m\oplus \gW+(\gS\cap C^{-1}\,\ima G_y)\oplus \gW$. We can write $\bsmat \x \\[1mm] \w\esmat=\bsmat \x_1 \\[1mm] \w_1\esmat+\bsmat \x_2 \\[1mm] \w_2\esmat$
with $\x_1\in \gV_m$, $\x_2\in \gS\cap C^{-1}\,\ima G_y$, so that $\x_2\in \gS$ and there exists $\xi$ such that $C\,\x_2=G_y\,\xi$. Now,
$\bsmat \x \\[1mm] \w \esmat={\bsmat \x_1\\[1mm] \w+\xi \esmat}+{\bsmat \x_2\\[1mm] -\xi\esmat}$, where $\bsmat \x_1\\[1mm] \w+\xi \esmat\in \gV_m\oplus \gW$ and $\bsmat \x_2\\[1mm] -\xi\esmat\in \gS\oplus \gW\,\cap\, \ker {[\,C\;\;G_y\,]}$.
Using (\ref{equltima}) 
\bea
(\gS+\gV_m)\oplus\gW \ns&\ns = \ns&\ns \left[
\left(\gS\cap \gV^\star_{\scriptscriptstyle \left(A,H,\bsmat \scriptscriptstyle{C} \\ \scriptscriptstyle{E} \esmat,\bsmat \scriptscriptstyle{G_y} \\[-0.5mm] \scriptscriptstyle{G_z} \esmat\right)}\right)+\gV_m \right]\oplus \gW \nonumber  \\
 \ns&\ns \subseteq \ns&\ns\left[
\left(\gS\cap C^{-1}\,\ima G_y \right)+\gV_m \right]\oplus \gW \nonumber  \\
 \ns&\ns= \ns&\ns \gV_m\oplus \gW+(\gS\cap C^{-1}\,\ima G_y)\oplus \gW,
\eea
since $\gV^\star_{\scriptscriptstyle \left(A,H,\bsmat \scriptscriptstyle{C} \\ \scriptscriptstyle{E} \esmat,\bsmat \scriptscriptstyle{G_y} \\[-0.5mm] \scriptscriptstyle{G_z} \esmat\right)}\subseteq C^{-1}\,\ima G_y$.
We have proved that $\gS+\gV_m$ is $(A,B,E,D_z)$-output nulling. Since $\gV_m$ is also $(A,B,E,D_z)$-self bounded, we have $\gV_m\supseteq \gV^\star_{\scriptscriptstyle (A,B,E,D_z)}\cap B\,\ker D_z$, so that also $\gS+\gV_m\supseteq \gV^\star_{\scriptscriptstyle (A,B,E,D_z)}\cap B\,\ker D_z$, i.e., $\gS+\gV_m$ is $(A,B,E,D_z)$-self bounded.
\endproof

Under the assumptions of the previous results, every $(A,H,C,G_y)$-input containing subspace $\gS$  such that $\gV_m\cap \gS_M\subseteq \gS\subseteq \gS_M$ is also $(A,H,C,G_y)$-self hidden. Indeed, from the obvious inclusion
$\gS_M=\gQ^\star_{\scriptscriptstyle \left(A,H,\bsmat \scriptscriptstyle{C} \\ \scriptscriptstyle{E} \esmat,\bsmat \scriptscriptstyle{G_y} \\[-0.3mm] \scriptscriptstyle{G_z} \esmat\right)}\subseteq \gQ^\star_{\scriptscriptstyle (A,H,C,G_y)}$,
we also have $\gS\subseteq \gS_M\subseteq \gQ^\star_{\scriptscriptstyle (A,H,C,G_y)}$.

We now seek to generalize the result in  
\cite[Lemma A.2]{DMCuellar-M-01}, which suggests that we can write, using the notation of this paper, $\gS+\gV_m
=\gR^\star_{\scriptscriptstyle (A,[\,B\;\;S\,],E,[\,D_z\;\;0\,])}$, where $S$ is basis matrix for $\gS$. However, it is easy to realize that in the biproper case this is not true. 
\begin{example}
{Consider e.g. $A=\bsmat 0 && 0 \\[1mm] 0 && -1 \esmat$, 
$B=\bsmat 1\\[1mm] 0\esmat$, $H=\bsmat 1\\[1mm] 1\esmat$, 
$C=[\,1\;\;0\,]$, $D_y=0$, $E=[\,0\;\;1\,]$, $D_z=1$ and $G_y=G_z=-1$. In this case $\gS_M=0_{\scriptscriptstyle \gX}$ and $\gV_m=\real^2$. Taking for example $\gS=\gS_M$ yields $\gS+\gV_m=\real^2$; however, a simple calculation shows that $\gR^\star_{\scriptscriptstyle (A,[\,B\;\;S\,],E,[\,D_z\;\;0\,])}=0_{\scriptscriptstyle \gX}$.}
\end{example}
The correct way to extend \cite[Lemma A.2]{DMCuellar-M-01} in the biproper case is to consider the subspace $\gR^\star_{\scriptscriptstyle (A,[\,B\;\;H \;\; S\,],E,[\,D_z\;\;G_z\;\; 0\,])}$ in place of  $\gR^\star_{\scriptscriptstyle (A,[\,B\;\;S\,],E,[\,D_z\;\;0\,])}$.  Correspondingly, in the case of zero feedthrough matrices one can show, by using the fact that (\ref{condalpha}-\ref{condbeta}) and $\gS^\star_{\scriptscriptstyle (A,H,C,G_y)} \subseteq \gV^\star_{\scriptscriptstyle (A,B,E,D_z)}$ are reduced to the chain of inclusions
$\ima H \subseteq \gS^\star_{\scriptscriptstyle (A,H,C,0)} \subseteq \gV^\star_{\scriptscriptstyle (A,B,E,0)}\subseteq \ker E$, one can easily see that 
$\gR^\star_{\scriptscriptstyle (A,[\,B\;\;H \;\; S\,],E,[\,0\;\;0\;\; 0\,])}=\gR^\star_{\scriptscriptstyle (A,[\,B\;\;S\,],E,[\,0\;\;0\,])}$.

 \begin{property}
 \label{lemA2m}
 Let (\ref{condalpha}-\ref{condbeta}) hold.
Let $\gS$ be an $(A,H,C,G_y)$-input containing subspace such that $\gV_m\cap \gS_M \subseteq \gS \subseteq \gS_M$, and let $S$ be a basis matrix for $\gS$. Then, $\gS^\star_{\scriptscriptstyle (A,H,C,G_y)} \subseteq \gV^\star_{\scriptscriptstyle (A,B,E,D_z)}$ if and only if
$\gS+\gV_m=\gR^\star_{\scriptscriptstyle (A,[\,B\;\;H \;\; S\,],E,[\,D_z\;\;G_z\;\; 0\,])}$.
\end{property}
\proof {\bf (Only if)}. We assume $\gS^\star_{\scriptscriptstyle (A,H,C,G_y)} \subseteq \gV^\star_{\scriptscriptstyle (A,B,E,D_z)}$. 
From Corollary~\ref{cor8}, $\gS_M\subseteq \gV^\star_{\scriptscriptstyle (A,B,E,D_z)}$. Since $\gS\subseteq \gS_M$, then also $\gS\subseteq \gV^\star_{\scriptscriptstyle (A,B,E,D_z)}$.
From Theorem~\ref{the12},  since $\gS\subseteq \gV^\star_{\scriptscriptstyle (A,B,E,D_z)}$, then $\gS\subseteq \gR^\star_{\scriptscriptstyle (A,[\,B\;\;S\,],E,[D_z\;\;0\,])}$, which is in turn contained in $\gR^\star_{\scriptscriptstyle (A,[\,B\;\;H \;\; S\,],E,[\,D_z\;\;G_z\;\; 0\,])}$. On the other hand, we have also $\gV_m=\gR^\star_{\scriptscriptstyle  A,[\,B\;\;H\,],E,[\,D_z\;\;G_z\,]}\subseteq \gR^\star_{\scriptscriptstyle (A,[\,B\;\;H \;\; S\,],E,[\,D_z\;\;G_z\;\; 0\,])}$, so that $\gS+\gV_m\subseteq \gR^\star_{\scriptscriptstyle  A,[\,B\;\;H\,],E,[\,D_z\;\;G_z\;\; 0\,]}$.
We know also that $\ima \bsmat H \\[1mm] G_z \esmat\subseteq \gV_m\oplus 0_{\scriptscriptstyle \gZ}+\ima \bsmat B \\[1mm] D_z \esmat$. Therefore, adding $\gS\oplus 0_{\scriptscriptstyle \gZ}$ on both sides yields
\bea
(\gS\oplus 0_{\scriptscriptstyle \gZ})+\ima \bmat{cc} H \\[-1.3mm] G_z \emat\subseteq \bigl((\gS+\gV_m) \oplus 0_{\scriptscriptstyle \gZ}\bigr)+\ima\bmat{cc} B \\[-1.3mm] D_z \emat.
\eea
In view of Property~\ref{pA4}, $\gS+\gV_m$ is $(A,B,E,D_z)$-self bounded, and it contains $\gS$. 
Since $\gS\subseteq  \gV^\star_{\scriptscriptstyle (A,B,E,D_z)}$ (recall that $\gS+\gV_m\in \Phi_{\scriptscriptstyle (A,B,E,D_z)}$), and since $\ima \bsmat H \\[1mm] G_z \esmat\subseteq (\gV^\star_{\scriptscriptstyle (A,B,E,D_z)}\oplus 0_{\scriptscriptstyle \gZ})+\ima\bsmat B \\[1mm] D_z \esmat$, then
$(\gS\oplus 0_{\scriptscriptstyle \gZ})+\ima \bsmat H \\[1mm] G_z \esmat\subseteq \bigl(\gV^\star_{\scriptscriptstyle (A,B,E,D_z)}\oplus 0_{\scriptscriptstyle \gZ}\bigr)+\ima\bsmat B \\[1mm] D_z \esmat$, 
so that by Theorem~\ref{prop3}  $\gR^\star_{\scriptscriptstyle (A,[\,B\;\;H \;\; S\,],E,[\,D_z\;\;G_z\;\; 0\,])}$ is the smallest $(A,B,E,D_z)$-self bounded subspace  such that 
\[
(\gS\oplus 0_{\scriptscriptstyle \gZ})+\ima \bmat{cc} H \\[-1.3mm] G_z \emat\subseteq \bigl(\gR^\star_{\scriptscriptstyle (A,[\,B\;\;H \;\; S\,],E,[\,D_z\;\;G_z\;\; 0\,])}\oplus 0_{\scriptscriptstyle \gZ}\bigr)+\ima\bmat{cc} B \\[-1.3mm] D_z \emat.
\]
Since it is the smallest,  $\gR^\star_{\scriptscriptstyle (A,[\,B\;\;H \;\; S\,],E,[\,D_z\;\;G_z\;\; 0\,])}\subseteq \gS+\gV_m$.\\ 
{\bf (If)}. It follows from $\gS^\star_{\scriptscriptstyle (A,H,C,G_y)}\subseteq \gS\subseteq \gS+\gV_m\subseteq \gV^\star_{\scriptscriptstyle (A,B,E,D_z)}$.
\endproof

We also present the dual of Property~\ref{lemA2m}.

 \begin{property}
 \label{A2dual}
 Let (\ref{condalpha}-\ref{condbeta}) hold.
Let $\gV$ be an  $(A,B,E,D_z)$-output nulling subspace such that $\gV_m \subseteq \gV \subseteq \gV_m+\gS_M$, and let $T$ be a full row-rank matrix such that $\ker T=\gV$. Then, $\gS^\star_{\scriptscriptstyle (A,H,C,G_y)} \subseteq \gV^\star_{\scriptscriptstyle (A,B,E,D_z)}$ if and only if
$\gV\cap \gS_M=\gQ^\star_{\scriptscriptstyle \left(A,H,\bsmat _{\scriptscriptstyle C}\\ _{\scriptscriptstyle E}\\ _{\scriptscriptstyle T} \esmat,\bsmat _{\scriptscriptstyle G_y} \\ _{\scriptscriptstyle G_z} \\ _{\scriptscriptstyle 0} \esmat\right)}$.
\end{property}

The previous result implies that $\gV\cap \gS_M\in \Psi_{\scriptscriptstyle (A,H,C,G_y)}$.
The following lemma, which was proved in \cite{DMCuellar-M-01} and \cite{DMCuellar-97} in the case of zero feedthrough matrices, relies 
  on a change of coordinate that was introduced in \cite[Theorem 5.2.2]{Basile-M-92}. However, here the structure of the system matrices in the new basis cannot be simplified as in \cite{Basile-M-92}. The following result shows how the line of attack of \cite[Theorem 5.2.2]{Basile-M-92} has to be changed substantially in the biproper case.

\begin{lemma}
\label{A3M}
 Let (\ref{condalpha}-\ref{condbeta}) hold.
Let $\gS^\star_{\scriptscriptstyle (A,H,C,G_y)} \subseteq \gV^\star_{\scriptscriptstyle (A,B,E,D_z)}$. 
Let $\gS$ be an $(A,H,C,G_y)$-input containing subspace such that $\gV_m\cap \gS_M \subseteq \gS \subseteq \gS_M$. 
Then
\bea
\sigma\Big(A+B\,F\,\Big|\,\frac{\gS+\gV_m}{\gV_m}\Big)=\sigma\Big(A+G\,C\,\Big|\,\frac{\gS}{\gV_m\cap\gS_M}\Big),\label{BMP}
\eea
where $F\in \mathfrak{F}_{\scriptscriptstyle (A,B,E,D_z)}(\gS+\gV_m)$ and $G\in \mathfrak{G}_{\scriptscriptstyle (A,H,C,G_y)}(\gV_m\cap \gS_M)$.
\end{lemma}

 \proof First, $\gS+\gV_m$ and $\gV_m$ are $(A,B,E,D_z)$-self bounded, see Property~\ref{pA4}, so that $\sigma\left(A+B\,F\,|\,\frac{\gS+\gV_m}{\gV_m}\right)$ is well defined for $F\in \mathfrak{F}_{\scriptscriptstyle (A,B,E,D_z)}(\gS+\gV_m)$. From Lemma~\ref{lem22}, $\gV_m\cap \gS_M$ is $(A,H,C,G_y)$-self hidden, so that $\sigma\left(A+G\,C\,|\,\frac{\gS}{\gV_m+\gS_M}\right)$ is well defined for $G\in \mathfrak{G}_{\scriptscriptstyle (A,H,C,G_y)}(\gV_m\cap \gS_M)$. 
Let $T=[\begin{array}{cccc} T_1 & T_2 & T_3 & T_4 \end{array}]$ be a nonsingular matrix such that
$\ima T_1 = \gV_m\cap \gS_M =\gQ^\star_{\scriptscriptstyle \left(A,H,\bsmat {\scriptscriptstyle C} \\ {\scriptscriptstyle E} \\ {\scriptscriptstyle T_m} \esmat,\bsmat {\scriptscriptstyle G_y} \\ {\scriptscriptstyle G_z} \\ {\scriptscriptstyle 0} \esmat\right)}$, $\ima [\begin{array}{cccc} T_1 & T_2  \end{array}] = \gV_m$, $\ima [\begin{array}{cccc} T_1 & T_3 \end{array}] = \gS$, $\ima [\begin{array}{cccc} T_1 & T_2 & T_3 \end{array}] = \gS+\gV_m=\gR^\star_{\scriptscriptstyle (A,[\,B\;\;H\;\;S\,],E,[\,D_z\;\;G_z\;\; 0\,])}$,
where $T_m$ is full row-rank and $\ker T_m=\gV_m$.\footnote{Notice that one or more of the matrices $T_i$ in $T$ may be empty.}
Let $\gV=\gS+\gV_m$. From Property~\ref{pA4}, $\gV=\gS+\gV_m\in \Phi_{\scriptscriptstyle (A,B,E,D_z)}$, and from Property~\ref{A2dual}, 
\bea
\begin{array}{rcl}
\gQ^\star_{\scriptscriptstyle \left(A,H,\bsmat {\scriptscriptstyle C} \\ {\scriptscriptstyle E} \\ {\scriptscriptstyle T_{\gV}} \esmat,\bsmat {\scriptscriptstyle G_y} \\ {\scriptscriptstyle G_z} \\ {\scriptscriptstyle 0} \esmat\right)} 
\ns&\ns = \ns&\ns \gV\cap \gS_M = (\gV_m\cap \gS_M)+\gS=\gS,
\end{array}
\eea
since $\gV_m\cap \gS_M\subseteq \gS$, 
where $T_{\gV}$ is full row-rank and $\ker T_{\gV}=\gV$.
Thus, $\ima [\begin{array}{cccc} T_1 & T_3 \end{array}] =\gQ^\star_{\scriptscriptstyle \left(A,H,\bsmat {\scriptscriptstyle C} \\ {\scriptscriptstyle E} \\ {\scriptscriptstyle T_{\gV}} \esmat,\bsmat {\scriptscriptstyle G_y} \\ {\scriptscriptstyle G_z} \\ {\scriptscriptstyle 0} \esmat\right)}=\gS$, from which we obtain
$\gS=\gV^\star_{\scriptscriptstyle \left(A,H,\bsmat {\scriptscriptstyle C} \\ {\scriptscriptstyle E} \\ {\scriptscriptstyle T_{\scriptscriptstyle \gV}} \esmat,\bsmat {\scriptscriptstyle G_y} \\ {\scriptscriptstyle G_z} \\ {\scriptscriptstyle 0} \esmat\right)}+\gS^\star_{\scriptscriptstyle \left(A,H,\bsmat {\scriptscriptstyle C} \\ {\scriptscriptstyle E} \\ {\scriptscriptstyle T_{\scriptscriptstyle  \gV}} \esmat,\bsmat {\scriptscriptstyle G_y} \\ {\scriptscriptstyle G_z} \\ {\scriptscriptstyle 0} \esmat\right)}$.
Since (\ref{condalpha}-\ref{condbeta}) and $\gS^\star_{\scriptscriptstyle (A,H,C,G_y)} \subseteq \gV^\star_{\scriptscriptstyle (A,B,E,D_z)}$ are assumed to hold, we have $\gS^\star_{\scriptscriptstyle \left(A,H,\bsmat \scriptscriptstyle{C} \\ \scriptscriptstyle{E} \esmat,\bsmat \scriptscriptstyle{G_y} \\[-0.5mm] \scriptscriptstyle{G_z} \esmat\right)}= \gS^\star_{\scriptscriptstyle \left(A,H,C,G_y\right)}\subseteq \gV$, so that from Theorem~\ref{the11d} we obtain
$\gS^\star_{\scriptscriptstyle \left(A,H,\bsmat {\scriptscriptstyle C} \\ {\scriptscriptstyle E} \\ {\scriptscriptstyle T_{\scriptscriptstyle \gV}} \esmat,\bsmat {\scriptscriptstyle G_y} \\ {\scriptscriptstyle G_z} \\ {\scriptscriptstyle 0} \esmat\right)}= \gS^\star_{\scriptscriptstyle \left(A,H,C,G_y\right)}$.
Clearly 
$
\gV^\star_{\scriptscriptstyle \left(A,H,\bsmat {\scriptscriptstyle C} \\ {\scriptscriptstyle E} \\ {\scriptscriptstyle T_{\scriptscriptstyle \gV}} \esmat,\bsmat {\scriptscriptstyle G_y} \\ {\scriptscriptstyle G_z} \\ {\scriptscriptstyle 0} \esmat\right)}\subseteq\bsmat C\\[1mm] E\\[1mm] T_{\gV}\esmat^{-1}\ima\bsmat G_y\\[1mm] G_z\\[1mm] 0\esmat\subseteq C^{-1}\ima G_y,
$
so that $\gS\subseteq C^{-1}\ima G_y +  \gS^\star_{\scriptscriptstyle \left(A,H,C,G_y\right)}$.
Now considering that
$
\gV_m\cap \gS_M =\gQ^\star_{\scriptscriptstyle \left(A,H,\bsmat {\scriptscriptstyle C} \\ {\scriptscriptstyle E} \\ {\scriptscriptstyle T_m} \esmat,\bsmat {\scriptscriptstyle G_y} \\ {\scriptscriptstyle G_z} \\ {\scriptscriptstyle 0} \esmat\right)}= \gV^\star_{\scriptscriptstyle \left(A,H,\bsmat {\scriptscriptstyle C} \\ {\scriptscriptstyle E} \\ {\scriptscriptstyle T_m} \esmat,\bsmat {\scriptscriptstyle G_y} \\ {\scriptscriptstyle G_z} \\ {\scriptscriptstyle 0} \esmat\right)}+\gS^\star_{\scriptscriptstyle \left(A,H,\bsmat {\scriptscriptstyle C} \\ {\scriptscriptstyle E} \\ {\scriptscriptstyle T_m} \esmat,\bsmat {\scriptscriptstyle G_y} \\ {\scriptscriptstyle G_z} \\ {\scriptscriptstyle 0} \esmat\right)}
$, adding $C^{-1}\ima G_y$ to both sides of the previous equation  and considering that by Theorem~\ref{the11d} there hold
$\gS^\star_{\scriptscriptstyle \left(A,H,\bsmat {\scriptscriptstyle C} \\ {\scriptscriptstyle E} \\ {\scriptscriptstyle T_m} \esmat,\bsmat {\scriptscriptstyle G_y} \\ {\scriptscriptstyle G_z} \\ {\scriptscriptstyle 0} \esmat\right)}=\gS^\star_{\scriptscriptstyle \left(A,H,C,G_y\right)}$
and
$
\gV^\star_{\scriptscriptstyle \left(A,H,\bsmat {\scriptscriptstyle C} \\ {\scriptscriptstyle E} \\ {\scriptscriptstyle T_m} \esmat,\bsmat {\scriptscriptstyle G_y} \\ {\scriptscriptstyle G_z} \\ {\scriptscriptstyle 0} \esmat\right)}\subseteq\bsmat C\\[1mm] E\\[1mm] T_{m}\esmat^{-1}\ima\bsmat G_y\\[1mm] G_z\\[1mm] 0\esmat\subseteq C^{-1}\ima G_y
$
we immediately obtain 
$C^{-1}\ima G_y+(\gV_m\cap \gS_M) = C^{-1}\ima G_y+\gS^\star_{\scriptscriptstyle \left(A,H,C,G_y\right)}$.
Finally 
\beann
\ima [\begin{array}{cccc} T_1 & T_3 \end{array}] \ns&\ns = \ns&\ns \gS\subseteq C^{-1}\ima G_y +  \gS^\star_{\scriptscriptstyle \left(A,H,C,G_y\right)}\\
\ns&\ns = \ns&\ns C^{-1}\ima G_y+(\gV_m\cap \gS_M)= C^{-1}\ima G_y+\ima T_1,
\eeann
which means that it is always possible to choose $T_3$ so that $\ima T_3\subseteq C^{-1}\ima G_y$.
If we denote by $\gR_{\scriptscriptstyle \gV_m+\gS}$ the output nulling reachability subspace on $\gV_m+\gS$, then $(\gS+\gV_m)\cap B\,\ker D_z=
\gR_{\scriptscriptstyle \gV_m+\gS}\cap B\,\ker D_z$. The subspace $\gR_{\scriptscriptstyle \gV_m+\gS}$ is contained $\gR^\star_{\scriptscriptstyle (A,B,E,D_z)}$, which 
in turn is contained in $\gV_m$. Thus,
$\ima [\begin{array}{ccc} T_1 & T_2 & T_3 \end{array}]\cap B\,\ker D_z = \ima [\begin{array}{ccc} T_1 & T_2\end{array}]\cap B\,\ker D_z$.
We can also choose $T_4$ in such a way that $\ima [\begin{array}{ccc} T_1 & T_2 & T_4 \end{array}]\supseteq B\,\ker D_z$.
 Let $A_1=T^{-1}\,A\,T$, $B_1=T^{-1} \,B$, $H_1=T^{-1}\,H$, $C_1=C\,T$, $E_1=E\,T$.
Let $\Omega= [\begin{array}{cc}\Omega_1 & \Omega_2 \end{array}]$ be a change of coordinate matrix in $\gU$ such that $\Omega_1$ is a basis for $B_1^{-1} (\gV_m+\gS)\cap \ker D_z$. We partition
$B_{\Omega} \defi B_1\,\Omega=\left[B_{i,j}\right]_{\tiny \substack{i = 1,\ldots,4 \\ j=1,2}}$
conformably with the change of basis $T$.
Since $B\,\ker D_z\subseteq \ima [\begin{array}{ccc} T_1 & T_2 & T_4\end{array}]$, we also have in the new basis
$B_1\,\ker D_z\subseteq \ima \bsmat I & 0 & 0 \\ 0 & I & 0 \\ 0 & 0 & 0\\ 0 & 0 & I \esmat$,
so that 
\[
\ima \bmat{cc} 
B_{11} \\[-1.3mm]
B_{21}\\[-1.3mm]
B_{31} \\[-1.3mm]
B_{41} \emat \subseteq  B_1\,\ker D_z \subseteq \ima \bmat{ccc} I & 0 & 0 \\[-1.3mm] 0 & I & 0 \\[-1.3mm] 0 & 0 & 0\\[-1.3mm] 0 & 0 & I \emat,
\]
which implies $B_{31}=0$. Partitioning
$A_1=T^{-1}A\,T=
\left[A_{i,j}\right]_{\tiny \substack{i= 1,\ldots,4 \\ j =1,\ldots,4}}$
and $F_{\Omega}=\Omega^{-1}\,\underbrace{F\,T}_{F_1}=
\left[F_{i,j}\right]_{\tiny \substack{i= 1,2 \\ j = 1,\ldots,4}}$
allows to partition the closed-loop matrix as
\beann
A_1\!+\!B_1 F_1=A_1\!+\!B_{\Omega} F_{\Omega}=\!
\left[ A_{i,j}\!+\!B_{i,1} F_{1,j}\!+\!B_{i,2} F_{2,j}\right]_{i,j=1,\ldots,4}.
\eeann
From the fact that $\gV_m$ is $(A+B\,F)$-invariant (in view of the self boundedness of $\gV_m$) 
we obtain
\beann
\begin{array}{rclccccrcl}
A_{31}\!+\!B_{31} F_{11}\!+\!B_{32} F_{21}\!\! \ns&\ns\! =\!\ns&\ns 0, &&&& A_{32}\!+\!B_{31} F_{12}\!+\!B_{32} F_{22}\!\!\ns&\ns\! =\!\ns&\ns 0, \\ 
A_{41}\!+\!B_{41} F_{11}\!+\!B_{42} F_{21}\!\!\ns&\ns\! =\!\ns&\ns 0,  &&&& A_{42}\!+\!B_{41} F_{12}\!+\!B_{42} F_{22} \!\! \ns&\ns\! =\!\ns&\ns 0.
\end{array}
\eeann
Moreover, since $\gV_m+\gS_M$ is $(A+B\,F)$-invariant, we also have
$A_{43}+B_{41}\,F_{13}+B_{42}\,F_{23}  =0$.
Thus, defining $A^F_{i,j}\defi A_{i,j}+B_{i,1}\,F_{1,j}+B_{i,2}\,F_{2,j}$ for $i,j \in \{1,\ldots,4\}$, we have
\[
A_1+B_1\,F_1=
\bmat{cccc}
A^F_{11}  & A^F_{12}  & A^F_{13} & A^F_{14} \\[-1.3mm]

A^F_{21} & A^F_{22}  & A^F_{23}  & A^F_{24}  \\[-1.3mm]

0 & 0& A_{33}+B_{32}\,F_{23}  & A_{34}+B_{32}\,F_{24}  \\[-1.3mm]

0  & 0  & 0  & A^F_{44} \emat,
\]
where we have used the fact that $B_{31}=0$. We now prove by contradiction that $B_{32}\,F_{23}=0$.
Recall that $\gS+\gV_m$ is $(A,B,E,D_z)$-output nulling, so that its reachability subspace is contained in $\gR^\star_{\scriptscriptstyle (A,B,E,D_z)}$; on the other hand, $\gV_m=\gR^\star_{\scriptscriptstyle (A,[\,B\;\;H\,],E,[\,D_z\;\;G_z\,])}\supseteq \gR^\star_{\scriptscriptstyle (A,B,E,D_z)}$. Hence, the reachability subspace on $\gS+\gV_m$ lies in $\gV_m$. 
Thus, the eigenvalues of $A+B\,F\,|\,\frac{\gS+\gV_m}{\gV_m}$ cannot be freely assigned by $F$. We have
\beann
\begin{array}{rcl}
\sigma(A\!+\!B F\,|\,\gS\!+\!\gV_m) \ns&\ns \!=\! \ns&\ns
\sigma(A\!+\!B F\,|\,\gV_m)\!\uplus  \sigma(A\!+\!B F\,|\,\frac{\gS+\gV_m}{\gV_m}) \\
 \ns&\ns \!=\! \ns&\ns \sigma(A\!+\!B F\,|\,\gV_m)\!\uplus \sigma (A_{33}\!+\!B_{32} F_{23}).
 \end{array}
\eeann
The submatrix $F_{23}$ in $A_1+B_1\,F_1$ does not affect $A+B\,F\,|\,\gV_m$. Thus, we can arbitrarily change $F_{23}$ without changing $\sigma(A+B\,F\,|\,\gV_m)$. In doing so, if $B_{32}\,F_{23}\neq 0$, we could change the eigenvalues of  $\sigma\left(A+B\,F\,|\,\frac{\gS+\gV_m}{\gV_m}\right)$ without changing $\sigma(A+B\,F\,|\,\gV_m)$. 
There holds
\beann
\begin{array}{rcl}
\sigma(A\!+\!B F\,|\,\gS\!+\!\gV_m) \ns&\ns = \ns&\ns\sigma(A\!+\!B F\,|\,\gR_{\scriptscriptstyle \gS+\gV_m})\\
\ns&\ns\ns&\ns \hspace{-1cm} \uplus  \sigma\left(A\!+\!B F\,\Big|\,\frac{\gV_m}{\gR_{\scriptscriptstyle \gS+\gV_m}}\right)\!\uplus \sigma\left(A\!+\!B F\,\Big|\,\frac{\gS+\gV_m}{\gV_m}\right)\!,
 \end{array}
\eeann
where the last two multisets are parts of the unreachable spectrum on $\gS+\gV_m$. If  we could change the eigenvalues of $\sigma\left(A+B\,F\,|\,\frac{\gS+\gV_m}{\gV_m}\right)$ without modifying those of $\sigma\left(A+B\,F\,\Big|\,\frac{\gV_m}{\gR_{\scriptscriptstyle \gS+\gV_m}}\right)$,  then we could change eigenvalues that are unreachable on $\gS+\gV_m$, leading to a contradiction.

Consider the following change of basis in the output space $\Upsilon=[\begin{array}{cc} \Upsilon_1 & \Upsilon_2 \end{array}]$, where $\Upsilon_1$ is a basis of $C_1\,\gS+\ima G_y$. 
We partition
$C_{\Upsilon}\defi \Upsilon^{-1}\,C_1=
\left[C_{i,j}\right]_{\tiny \substack{i=1,2 \\ j=1,\ldots,4}}$.
Since $\ima T_3\subseteq C^{-1}\,\ima G_y$, in the new basis we have
$\ima \bsmat 0\\[0.5mm] 0\\[0.5mm] I\\[0.5mm] 0 \esmat \subseteq C_1^{-1}\,\ima G_y$,
so that 
$
\ima \left(C_1\,\bsmat 0\\[0.5mm] 0\\[0.5mm] I\\[0.5mm] 0 \esmat\right) \subseteq \ima G_y$.
Pre-multiplying by $\Upsilon^{-1}$ we obtain
$\bsmat C_{13}\\[1mm] C_{23} \esmat=\Upsilon^{-1}\ima G_y=\bsmat \star \\[1mm] 0 \esmat$,
since $\ima G_y\subseteq \ima \Upsilon_1$. 
Thus, $C_{23}=0$.
Moreover, we partition $G$ conformably as
$G_{\Upsilon}=\underbrace{T^{-1}\,G}_{G_1}\,\Upsilon=
\left[G_{i,j}\right]_{\tiny \substack{i = 1,\ldots,4 \\ j = 1,2} }$.
Thus, defining $A^G_{i,j}\defi A_{i,j}+G_{i,1}\,C_{1,j}+G_{i,2}\,C_{2,j}$, we can write $A_G \defi A_1+G_1\,C_1=A_1+G_{\Upsilon}\,C_{\Upsilon}=
\left[A^G_{i,j}\right]_{\tiny \substack{i= 1,\ldots,4 \\ j =1,\ldots,4} }$.
Since $\gS$ and $\gS_M\cap \gV_m$ are $(A,H,C,G_y)$-self hidden, the input containing friend $G$ is also an $(A,H,C,G_y)$-input containing friend for $\gS$, and in this basis
$\ima \left(A_G\,\bsmat I \\[0mm] 0\\[0mm] 0\\[0mm] 0\esmat\right)\subseteq \ima \bsmat I \\[0mm] 0\\[0mm] 0\\[0mm] 0\esmat$
 and $\ima \left(A_G\,\bsmat I &0\\[0mm] 0&0\\[0mm] 0&I\\[0mm] 0&0\esmat\right)\subseteq \ima \bsmat I &0\\[0mm] 0&0\\[0mm] 0&I\\[0mm] 0&0\esmat$
give
\[
A_1+G_1\,C_1=A_1+G_{\Upsilon}\,C_{\Upsilon}=\bmat{cccc}
A^G_{11} & A^G_{12}  & A^G_{13} & A^G_{14} \\[-1.3mm]

0 & A^G_{22} & 0  & A^G_{24}\\[-1.3mm]

0 & A^G_{32} & A_{33}+G_{31}\,C_{13} & A^G_{34} \\[-1.3mm]

0 & A^G_{42} & 0& A^G_{44}  \emat\!,
\]
since, as shown above, $C_{23}=0$.
We have
\beann
\begin{array}{rcl}
&&\hspace{-3mm} \sigma(A_1\!+\!G_1 C_1) \!=\! \sigma(A_{11}\!+\!G_{11} C_{11}\!+\!G_{12} C_{21})\!\uplus \sigma(A_{33}\!+\!G_{31} C_{13}) \\
&& \; \uplus \, \sigma\left(\bmat{cccc}
A_{22}\!+\!G_{21} C_{12}\!+\!G_{12} C_{22} & A_{24}\!+\!G_{21} C_{14}\!+\!G_{22} C_{24} \\
 A_{42}\!+\!G_{41} C_{12}\!+\!G_{42} C_{22}  & A_{44}\!+\!G_{41} C_{14}\!+\!G_{42} C_{24}  \emat \right)\!\!.
 \end{array}
 \eeann
 We show that $G_{31}\,C_{13}=0$. To this end, notice that $\gS\supseteq \gS^\star_{\scriptscriptstyle (A,H,C,G_y)}$ and $\gS_M\cap \gV_m \supseteq \gS^\star_{\scriptscriptstyle (A,H,C,G_y)}$. This implies that the eigenvalues of $A+G\,C\,|\,\frac{\gS}{\gS_M\cap \gV_m}$ do not depend on $G$. Since 
$\sigma(A+G\,C\,|\,\frac{\gS}{\gS_M\cap \gV_m})=\sigma(A_{33}+G_{31}\,C_{13})$, and since $G_{31}$ does not operate on the eigenvalues of $A+G\,C\,|\,\gS_M\cap \gV_m$,
we have $G_{31}\,C_{13}=0$. 
From $\sigma\left(A+B\,F\,|\,\frac{\gS+\gV_m}{\gV_m}\right)=\sigma(A_{33})$ and $\sigma\left(A+G\,C\,|\,\frac{\gS}{\gV_m\cap\gS_M}\right)=\sigma(A_{33})$, we see that (\ref{BMP}) holds.
\endproof

The following result generalizes \cite[Lemma A.5]{DMCuellar-M-01} to the case of possibly non-zero feedthrough matrices. 

 \begin{property}
 \label{A5m}
 Let (\ref{condalpha}-\ref{condbeta}) hold. Let $\gS^\star_{\scriptscriptstyle (A,H,C,G_y)} \subseteq \gV^\star_{\scriptscriptstyle (A,B,E,D_z)}$. Then
\begin{enumerate}
\item $\langle A\,|\,\ima B +\ima H \rangle=\gR+\gV_m$;
\item $\langle \ker C \cap \ker E \,|\,A\rangle=\gQ \cap \gS_M$.
\end{enumerate}
\end{property}
\proof
We first prove that $\langle A\,|\,\ima B +\ima H \rangle\subseteq \langle A\,|\,\ima B \rangle+\gV_m$. We observe that 
\beann
\begin{array}{rcl}
\langle A\,|\,\ima B +\ima H \rangle \ns&\ns = \ns&\ns \langle A\,|\,\ima [\begin{array}{cc} B & H \end{array}] \rangle  \\
\ns&\ns = \ns&\ns \ima [\begin{array}{cc} B & H \end{array}] +\ima \bigl(A\, [\begin{array}{cc} B & H \end{array}]\bigr) + \ldots +\ima \bigl(A^{n-1}\, [\begin{array}{cc} B & H \end{array}]\bigr) \\
\ns&\ns = \ns&\ns \ima B+\ima H+\ima [\begin{array}{cc} A\,B & A\,H \end{array}]+\ldots+
\ima [\begin{array}{cc} A^{n-1}\,B & A^{n-1}\,H \end{array}]\\
\ns&\ns = \ns&\ns \ima B+\ima H+\ima (A\,B)+\ima (A\,H)+\ldots+\ima (A^{n-1}\,B)+\ima (A^{n-1}\,H)\\
\ns&\ns = \ns&\ns  \ima [\begin{array}{cccc} B & A\,B & \ldots & A^{n-1}\,B \end{array}]+\ima [\begin{array}{cccc} H & A\,H & \ldots & A^{n-1}\,H \end{array}] \\
\ns&\ns = \ns&\ns \langle A\,|\,\ima B \rangle+\langle A\,|\,\ima H \rangle.
\end{array}
\eeann
Recall that, since conditions (\ref{condalpha}-\ref{condbeta}) hold, then also from Theorem~\ref{the14} we have
$\ima \bsmat H \\[1mm] G_z \esmat \subseteq (\gV_m  \oplus 0_{\scriptscriptstyle \gZ})+\ima \bsmat B\\[1mm] D_z \esmat$.
This implies in particular that $\ima H\subseteq \gV_m+\ima B$, so that also
\beann
\begin{array}{rcl}
A^{n-1}\ima H \ns&\ns \subseteq \ns&\ns  A^{n-1}\,\gV_m+\ima (A^{n-1}\,B) \\
A^{n-2}\ima H \ns&\ns \subseteq \ns&\ns  A^{n-2}\,\gV_m+\ima (A^{n-2}\,B) \\
\ns&\ns \vdots \ns&\ns \\
A\ima H \ns&\ns \subseteq \ns&\ns  A\,\gV_m+\ima (A\,B).
\end{array} 
\eeann
Moreover, again from Theorem~\ref{the14}, $\gV_m$ is $(A,[\begin{array}{cc} B & H \end{array}],E,[\begin{array}{cc}D_z & 0 \end{array}])$-output nulling, so that, in particular, $\gV_m$ is $(A,[\begin{array}{cc} B & H \end{array}])$-controlled invariant, i.e,  
$A\,\gV_m\subseteq \gV_m+\ima [\begin{array}{cc} B & H \end{array}]$, which leads to 
\beann
\begin{array}{rcl}
A^n\,\gV_m \ns&\ns \subseteq \ns&\ns A^{n-1}\bigl(\gV_m+\ima B+\ima H\bigr) \\
A^{n-1}\,\gV_m \ns&\ns \subseteq \ns&\ns A^{n-2}\bigl(\gV_m+\ima B+\ima H\bigr) \\
\ns&\ns \vdots \ns&\ns \\
A\,\gV_m \ns&\ns \subseteq \ns&\ns \gV_m+\ima B+\ima H.
\end{array}
\eeann
We have 
\beann
\begin{array}{rcl}
\langle A\,|\,\ima B +\ima H \rangle \ns&\ns = \ns&\ns \langle A\,|\,\ima B \rangle+\langle A\,|\,\ima H \rangle \\
 \ns&\ns = \ns&\ns \langle A\,|\,\ima B \rangle+\ima H+\ima (A\,H)+\ldots+\ima (A^{n-1}\,H)\\
  \ns&\ns \subseteq  \ns&\ns \langle A\,|\,\ima B \rangle+\ima H+\ima (A\,H)+\ldots+\ima (A^{n-1}\,H)+A^n\,\gV_m\\
   \ns&\ns \subseteq  \ns&\ns \langle A\,|\,\ima B \rangle+\ima H+\ima (A\,H)+\ldots+\ima (A^{n-1}\,H)+A^{n-1}\,\bigl(\gV_m+\ima B+\ima H\bigr)\\      \ns&\ns \subseteq  \ns&\ns \langle A\,|\,\ima B \rangle+\ima H+\ima (A\,H)+\ldots+\ima (A^{n-1}\,H)+A^{n-1}\,\gV_m \\
      \ns&\ns \subseteq  \ns&\ns \langle A\,|\,\ima B \rangle+\ima H+\ima (A\,H)+\ldots+\ima (A^{n-2}\,H)+A^{n-1}\,\gV_m 
\\
   \ns&\ns \subseteq  \ns&\ns \langle A\,|\,\ima B \rangle+\ima H+\ima (A\,H)+\ldots+\ima (A^{n-2}\,H)+A^{n-2}\,\bigl(\gV_m+\ima B+\ima H\bigr)
\\
      \ns&\ns \subseteq  \ns&\ns \langle A\,|\,\ima B \rangle+\ima H+\ima (A\,H)+\ldots+\ima (A^{n-2}\,H)+A^{n-2}\,\gV_m 
            \end{array} \eeann
  \beann
\begin{array}{rcl}   
      \ns&\ns \subseteq  \ns&\ns \langle A\,|\,\ima B \rangle+\ima H+\ima (A\,H)+\ldots+\ima (A^{n-3}\,H)+A^{n-1}\,\gV_m 
\\
      \ns&\ns \vdots \ns&\ns \\
            \ns&\ns \subseteq  \ns&\ns 
       \langle A\,|\,\ima B \rangle+\gV_m.
       \end{array}
\eeann
We now need to prove the opposite inclusion, i.e., $\langle A\,|\,\ima B +\ima H \rangle\supseteq \langle A\,|\,\ima B \rangle+\gV_m$. Obviously $\langle A\,|\,\ima B \rangle \subseteq \langle A\,|\,\ima B +\ima H \rangle$ and $\gV_m \subseteq \langle A\,|\,\ima B +\ima H \rangle$ (recall that $\gV_m=\gR^\star_{\scriptscriptstyle (A,[\, B \;\; H \,],E,[\,E \;\; G_z \,])}$, so that in particular it is contained in $\langle A\,|\,\ima B +\ima H \rangle$). Thus, $\langle A\,|\,\ima B +\ima H \rangle\supseteq \langle A\,|\,\ima B \rangle+\gV_m$ holds.\\
The second statement can be proved by duality.
\endproof

\section{Fixed poles}
\label{sec:fixed}
Given an $(A,B,E,D_z)$-output nulling subspace $\gV$ and an $(A,H,C,G_y)$-input containing subspace $\gS$, we define
\bea
 \sigma_{\rm fixed}(\gS,\gV)\defi\sigma_{\rm fixed}(\gV)\uplus\sigma_{\rm fixed}(\gS).
\eea
The key idea behind the approach of \cite{DMCuellar-M-01}, on which our development also hinges, is to show that:
\begin{enumerate}
\item If DDPDOF is solvable and the resulting closed-loop matrix is $\widehat{A}$, it is also always solvable by using a solution triple $(\gS,\gV;K)$ where $\gV$ is $(A,B,E,D_z)$-self bounded and $\gS$ is $(A,H,C,G_y)$-self hidden which improve on the original solution, in the sense that the fixed poles associated with the triple $(\gS,\gV;K)$ are contained in the closed-loop eigenvalues $\sigma(\widehat{A})$;
\item We show that there exist ``minimal solution triples'' whose fixed poles are the same, and are always contained in the fixed poles associated with any other solution triple $(\gS,\gV;K)$ where $\gV$ is $(A,B,E,D_z)$-self bounded and $\gS$ is $(A,H,C,G_y)$-self hidden. 
\end{enumerate}

Thus, 
looking for the minimal solution triples in the set of solution triples $(\gS,\gV;K)$ where $\gV$ is $(A,B,E,D_z)$-self bounded and $\gS$ is $(A,H,C,G_y)$-self hidden does not cause any loss of generality.
Every controller that solves DDPDOF will give rise to a closed-loop eigenstructure which contains the fixed poles of the ``minimal solution triples'', even if such controller has not been constructed by using such minimal solution triples.
The following result has been proved for biproper systems in \cite[Lemma 4]{N-TAC-08}. 

\begin{lemma}
	\label{thm:morse}
Let $\gV$ and $\gV^\prime$ be two $(A,B,C,D)$-output nulling subspaces such that
$\gV=\gV^\prime+\gR_{\gV}$.
There exists $F\in \mathfrak{F}_{\scriptscriptstyle (A,B,C,D)}(\gV)$ such that $\sigma\left(A+B\,F\,|\,\frac{\gV}{\gV^\prime}\right)$ is freely assignable. Moreover, for all $F^\prime \in \mathfrak{F}_{\scriptscriptstyle (A,B,C,D)}(\gV^\prime)$ and for all $F \in \mathfrak{F}_{\scriptscriptstyle (A,B,C,D)}(\gV)$ we have
\bea
\sigma \Big(A+B\,F\,\Big|\,\frac{\gV}{\gR_{\gV}}\Big) \subseteq \sigma \Big(A+B\,F^\prime\,\Big|\,\frac{\gV^\prime}{\gR_{\gV^\prime}}\Big).
\eea
\end{lemma}
We are ready to present the following generalization of   \cite[Lemma 5]{DMCuellar-M-01}.  Here the situation is, however, substantially different. In fact, here we also need to prove the existence of a matrix $K$ that renders the closed-loop system well posed. For its proof, we recall the definition of 
projection and the intersection of subspaces in extended vector spaces. Let $\gT$ be a subspace of
$\gX\oplus \gP$. The projection of $\gT$ on $\gX$ is defined as $\mathfrak{p}\,(\gT) \defi \left\{\x\in \gX\,\Big|\,\exists \,\p  \in \gP\,:\; \bsmat \x \\[1mm] \p \esmat \in \gT\right\}$ and the intersection of $\gT$ with $\gX$ is defined as $\mathfrak{i}\,(\gT) \defi \left\{\x\in \gX\,\Big|\, \bsmat \x \\[1mm] {0} \esmat \in \gT\right\}$, see \cite{Basile-M-92,Trentelman-SH-01,PN-submitted} for details.

\begin{lemma}
\label{giapponesi}
Assume that DDPDOF is solvable.
 Let  $\sigma(\widehat{A})$ be the multiset of closed-loop eigenvalues. There exists a solution triple $(\gS,\gV;K)$ for DDPDOF such that $\gV$ is $(A,B,E,D_z)$-self bounded and $\gS$ is $(A,H,C,G_y)$-self hidden and such that $\sigma_{\rm fixed}(\gS,\gV)\subseteq\sigma(\widehat{A})$.
\end{lemma}

\proof
Since DDPDOF, is solvable, there exists an $\widehat{A}$-invariant subspace $\widehat{\gI}$ contained in $\ker \widehat{C}$ such that $\gS_{\mathfrak i}=\mathfrak{i}(\widehat{\gI})$ is $(A,H,C,G_y)$-input containing 
and $\gV_{\mathfrak p}=\mathfrak{p}(\widehat{\gI})$ is $(A,B,E,D_z)$-output nulling, $\gS_{\mathfrak i}\subseteq \gV_{\mathfrak p}$ (see \cite[Thm. 4.6]{Trentelman-SH-01} and  \cite{PN-submitted} as well as Appendix A), and satisfy 
\bea
\label{cond1}
 \ima \bmat{cc} H \\[-1.3mm] G_z \emat & \subseteq &  (\gV_{\mathfrak p} \oplus {0}_{\gZ}) +\ima \bmat{cc} B \\[-1.3mm]  D_z\emat;\\ 
\label{cond2}
\ker \,[\begin{array}{cc} E & G_z \end{array}]  & \supseteq &  (\gS_{\mathfrak i} \oplus \gW)  \cap \ker \,[\begin{array}{cc}  C & G_y \end{array}]. 
\eea
We need to prove that there exists $K$ such that (\ref{condalpha}-\ref{condbeta}) holds with $\gV_{\mathfrak p}$ and $\gS_{\mathfrak i}$ in place of $\gV^\star_{\scriptscriptstyle (A,B,E,D_z)}$ and $\gS^\star_{\scriptscriptstyle (A,H,C,G_y)}$.
Let $K=D_c\,W=D_c\,(I-D_y\,D_c)^{-1}$. The matrix inversion lemma ensures that $I+K\,D_y$ is invertible. 
It remains to prove that $K$ satisfies (\ref{clK}).
 Rewriting (\ref{clK}) using 
 $K=D_c\,(I-D_y\,D_c)^{-1}$ yields
\bea
\label{clKnew1}
 \bmat{cc} A\!+\!B D_c W C & H\!+\!B D_c W G_y \\[-1.3mm] 
 E\!+\!D_z D_c W C & G_z\!+\! D_z D_c W G_y \emat\!\!(\gS_{\mathfrak i}  \oplus \gW)\subseteq \gV_{\mathfrak p} \oplus 0_{\scriptscriptstyle \gZ}.
 \eea
 Let $\bsmat \s \\[1mm] \w \esmat \in \gS_{\mathfrak i}  \oplus \gW$. 
 We want to prove that 
 \bea
 \label{toprove}
  \bmat{cc} A\,\s+B\,D_c\,W\,C\,\s+ H\,\w+B\,D_c\,W\,G_y\,\w \\[-1.3mm] 
 E\,\s+D_z\,D_c\,W\,C\,\s+ G_z\,\w+ D_z\,D_c\,W\,G_y\,\w \emat\in \gV_{\mathfrak p} \oplus 0_{\scriptscriptstyle \gZ}.
 \eea
 Since $\s\in \gS_{\mathfrak i}=\mathfrak{i}(\widehat{\gI})$, we have $\bsmat \s \\[1mm] 0 \esmat\in \widehat{\gI}$, and from the $\widehat{A}$-invariance of $\widehat{\gI}$ we find
$\widehat{A}\,\bsmat \s \\[1mm] 0 \esmat=\bsmat
 A\,\s+B\,D_c\,W\,C\,\s \\[1mm] E\,\s+D_z\,D_c\,W\,C\,\s \esmat\in \widehat{\gI}$.
 It follows that  $A\,\s+B\,D_c\,W\,C\,\s \in \mathfrak{p}(\widehat{\gI})=\gV_{\mathfrak p}$. Moreover, since $\ima \widehat{H}\subseteq \widehat{\gI}$, we have 
$\bsmat H\,\w+B\,D_c\,W\,G_y\,\w \\[1mm] B_c\,W G_y \,\w \esmat \in \widehat{\gI}$,
so that, in particular, $H\,\w+B\,D_c\,W\,G_y\,\w\in \mathfrak{p}(\widehat{\gI})=\gV_{\mathfrak p}$. 
Since the system is disturbance decoupled, the feedthrough $G_z+ D_z\,D_c\,(I-D_y\,D_c)^{-1}\,G_y$ in (\ref{clKnew1}) is zero. Hence, it remains to show that $E\,\s+D_z\,D_c\,(I-D_y\,D_c)^{-1}\,C\,\s=0$. This follows from $\widehat{C}\,\widehat{\gI}=0_{\scriptscriptstyle \gZ}$, which gives $\widehat{C}\,\bsmat \s \\[1mm] 0 \esmat=0$. This yields immediately $E\,\s+D_z\,D_c\,(I-D_y\,D_c)^{-1}\,C\,\s=0$. We have proved (\ref{toprove}).\\
Let us define $\gV\defi\gV_{\mathfrak p}+\gR^\star_{\scriptscriptstyle (A,B,E,D_z)}$ and $\gS\defi \gS_{\mathfrak i}\cap\gQ^\star_{\scriptscriptstyle (A,H,C,G_y)}$.
Obviously, $\gS$ is $(A,H,C,G_y)$-input containing,
$\gV$ is $(A,B,E,D_z)$-output nulling, and $\gS\subseteq \gV$. Moreover (\ref{cond1}-\ref{cond2}) are satisfied. 
Furthermore, if $K$ is such that $\gS_{\mathfrak i}$, $\gV_{\mathfrak p}$ and $K$ form a solution triple $(\gS_{\mathfrak i},\gV_{\mathfrak p};K)$, the chain of inclusions 
\beann
&&\hspace{-3mm} \bmat{cc} \! A\!+\!B K C \!& H\!+\!B K G_y\! \\[-1.3mm] 
\! E\!+\!D_z K C \!& G_z\!+\! D_z K G_y\! \emat\!(\gS \oplus \gW)\\
&&\hspace{3mm}
\!\subseteq \! \bmat{cc} \! A\!+\!B K C \!& H\!+\!B K G_y\! \\[-1.3mm] 
\! E\!+\!D_z K C \!& G_z\!+\! D_z K G_y\!\emat\!\!(\gS_{\mathfrak i} \oplus \gW)\!\subseteq\! \gV_{\mathfrak p} \oplus 0_{\scriptscriptstyle \gZ}\!\subseteq\! \gV \! \oplus 0_{\scriptscriptstyle \gZ}
\eeann
guarantees that $(\gS,\gV;K)$ is a solution triple.
Finally, by construction, $\gV$ is $(A,B,E,D_z)$-self bounded because it contains $\gR^\star_{\scriptscriptstyle (A,B,E,D_z)}$, and 
 $\gS$ is $(A,H,C,G_y)$-self hidden because it is contained in $\gQ^\star_{\scriptscriptstyle (A,H,C,G_y)}$.
Now we show that $\sigma_{\rm fixed}(\gS,\gV)\subseteq\sigma(\widehat{A})$. We only need to prove that $\sigma_{\rm fixed}(\gV)\subseteq\sigma(\widehat{A})$,  because  $\sigma_{\rm fixed}(\gS)\subseteq\sigma(\widehat{A})$ follows by the self duality of the problem. \\
To prove that $\sigma_{\rm fixed}(\gV)\subseteq\sigma(\widehat{A})$, we first prove that $\sigma_{\rm fixed}(\gV) \subseteq \sigma_{\rm fixed}(\gV_{\mathfrak p})$ and then we prove that $\sigma_{\rm fixed}(\gV_{\mathfrak p})\subseteq\sigma(\widehat{A})$.
Using (\ref{sfixV}), for $F\in \mathfrak{F}_{\scriptscriptstyle (A,B,E,D_z)}(\gV)$ and $F_{\mathfrak p}\in \mathfrak{F}_{\scriptscriptstyle (A,B,E,D_z)}(\gV_{\mathfrak p})$
\beann
\sigma_{\rm fixed}(\gV) &=& \sigma \Big(A_F\,\Big|\,\frac{\gX}{\langle A\,|\,\ima B \rangle}\Big) \uplus 
\sigma \Big(A_F\,\Big|\,\frac{\gV\cap \langle A\,|\,\ima B \rangle}{\gR_{\gV}}\Big),\\
\sigma_{\rm fixed}(\gV_{\mathfrak p}) &=& \sigma \Big(A_{F_{\mathfrak p}}\,\Big|\,\frac{\gX}{\langle A\,|\,\ima B \rangle}\Big) \uplus 
\sigma \Big(A_{F_{\mathfrak p}}\,\Big|\,\frac{\gV_{\mathfrak p}\cap \langle A\,|\,\ima B \rangle}{\gR_{\gV_{\mathfrak p}}}\Big),
\eeann
where $A_F=A+B\,F$ and $A_{F_{\mathfrak p}}=A+B\,F_{\mathfrak p}$.
Since $\sigma \left(A_F\,\Big|\,\frac{\gX}{\langle A\,|\,\ima B \rangle}\right)$ does not depend on $F$ we have 
\bea
\sigma \Big(A_F\,\Big|\,\frac{\gX}{\langle A\,|\,\ima B \rangle}\Big)=\sigma \Big(A_{F_{\mathfrak p}}\,\Big|\,\frac{\gX}{\langle A\,|\,\ima B \rangle}\Big)\!,
\eea
so that it is enough to prove that for all $F\in \mathfrak{F}_{\scriptscriptstyle (A,B,E,D_z)}(\gV)$ and $F_{\mathfrak p}\in \mathfrak{F}_{\scriptscriptstyle (A,B,E,D_z)}(\gV_{\mathfrak p})$ 
\bea
\sigma\Big(\!A_F\,\Big|\,\frac{\gV\cap\langle A\,| \ima B \rangle}{\gR_{\gV}}\Big) \!\subseteq 
\sigma\Big(\!A_{F_{\mathfrak p}}\,\Big|\,\frac{\gV_{\mathfrak p}\cap\langle A\,|\ima B \rangle}{\gR_{\gV_{\mathfrak p}}}\Big).
\eea
From $\gR^\star_{\scriptscriptstyle (A,B,E,D_z)} \subseteq\langle A\,|\,\ima B \rangle$ and $\gV= \gV_{\mathfrak p}+\gR^\star_{\scriptscriptstyle (A,B,E,D_z)}$, we find
\bea
\gV\cap\langle A\,|\,\ima B \rangle \ns&\ns = \ns&\ns (\gV_{\mathfrak p}+\gR^\star_{\scriptscriptstyle (A,B,E,D_z)}) \cap \langle A\,|\,\ima B \rangle\nonumber \\
\ns&\ns = \ns&\ns (\gV_{\mathfrak p}\cap \langle A\,|\,\ima B \rangle) + \gR^\star_{\scriptscriptstyle (A,B,E,D_z)}. 
\eea
From $\gV= \gV_{\mathfrak p}+\gR^\star_{\scriptscriptstyle (A,B,E,D_z)}$ we find that $\gR_{\gV}$ contains $\gR^\star_{\scriptscriptstyle (A,B,E,D_z)}$. However, $\gR^\star_{\scriptscriptstyle (A,B,E,D_z)}$ is the largest $(A,B,E,D_z)$-reachability output nulling subspace, which implies $\gR^\star_{\scriptscriptstyle (A,B,E,D_z)}=\gR_{\gV}$.
From the obvious inclusion $\gR^\star_{\scriptscriptstyle (A,B,E,D_z)}\subseteq \langle A\,|\,\ima B\rangle$ the $(A,B,E,D_z)$-reachability subspace on $\gV\cap\langle A\,|\,\ima B \rangle$  
is exactly $\gR^\star_{\scriptscriptstyle (A,B,E,D_z)}$, i.e, $\gR^\star_{\scriptscriptstyle (A,B,E,D_z)}=\gR_{\gV\cap\langle A\,|\,\ima B \rangle}$, and therefore
$\gV\cap\langle A\,|\,\ima B \rangle=(\gV_{\mathfrak p}\cap\langle A\,|\,\ima B \rangle)+\gR_{\gV\cap\langle A\,|\,\ima B \rangle}$.
We can apply Lemma~\ref{thm:morse} to the previous subspaces and we find that for all $F\in \mathfrak{F}_{\scriptscriptstyle (A,B,E,D_z)}(\gV)$ and for all $F_{\mathfrak p}\in\mathfrak{F}_{\scriptscriptstyle (A,B,E,D_z)}(\gV_{\mathfrak p})$
\bea
\sigma\left(\!A_F\,\Big| \frac{\gV\cap\langle A\,|\ima B \rangle}{\gR_{\gV\cap\langle A\,|\,\ima B \rangle}}\right)\!
\subseteq\sigma\Big(\!A_{F_{\mathfrak p}}\,\Big|\, \frac{\gV_{\mathfrak p}\cap\langle A\,|\ima B \rangle}{\gR_{\gV_{\mathfrak p}\cap\langle A\,|\ima B \rangle}}\Big).
\eea
Since clearly $\gR_{\gV\cap\langle A\,|\,\ima B\rangle}=\gR_{\gV}$ and $\gR_{\gV_{\mathfrak p}\cap\langle A\,|\,\ima B \rangle}=\gR_{\gV_{\mathfrak p}}$, then for all $F\in \mathfrak{F}_{\scriptscriptstyle (A,B,E,D_z)}(\gV)$ and for all $F_{\mathfrak p}\in\mathfrak{F}_{\scriptscriptstyle (A,B,E,D_z)}(\gV_{\mathfrak p})$
\bea
\sigma\left(\!A_F\,\Big|\, \frac{\gV\cap\langle A\,|\ima B \rangle}{\gR_{\gV}}\right)\!\subseteq 
\sigma\Big(\!A_{F_{\mathfrak p}}\,\Big|\, \frac{\gV_{\mathfrak p}\cap\langle A\,|\,\ima B \rangle}{\gR_{\gV_{\mathfrak p}}}\Big),
\eea
so that $\sigma_{\rm fixed}(\gV) \subseteq \sigma_{\rm fixed}(\gV_{\mathfrak p})$.
Now we prove that $\sigma_{\rm fixed}(\gV_{\mathfrak p})\subseteq\sigma(\widehat{A})$. 
Using \eqref{sfixV} we have for $F_{\mathfrak p}\in \mathfrak{F}_{\scriptscriptstyle (A,B,E,D_z)}(\gV_{\mathfrak p})$
\bea
\sigma_{\rm fixed}(\gV_{\mathfrak p}) = \sigma\Big(A_{F_{\mathfrak p}}\,\Big|\,\frac{\gX}{\langle A\,|\,\ima B \rangle+\gV_{\mathfrak p}}\Big) \uplus 
\sigma \Big(A_{F_{\mathfrak p}}\,\Big|\,\frac{\gV_{\mathfrak p}}{\gR_{\gV_{\mathfrak p}}}\Big).
\eea
Clearly $\sigma\left(A_{F_{\mathfrak p}}\,\Big|\,\frac{\gX}{\langle A\,|\,\ima B \rangle+\gV_{\mathfrak p}}\right) \subseteq \sigma\left(A_{F_{\mathfrak p}}\,\Big|\,\frac{\gX}{\langle A\,|\,\ima B \rangle}\right)\subseteq \sigma(\widehat A)$, and thus it is enough to prove that $\sigma\left(A_{F_{\mathfrak p}}\Big|\frac{\gV_{\mathfrak p}}{\gR_{\gV_{\mathfrak p}}}\right)\!\subseteq\! \sigma(\widehat A)$ for all  $F_{\mathfrak p}\in \mathfrak{F}_{\scriptscriptstyle (A,B,E,D_z)}(\gV_{\mathfrak p})$. 
In particular, we show that $\sigma\left(A_{F_{\mathfrak p}}\Big|\frac{\gV_{\mathfrak p}}{\gR_{\gV_{\mathfrak p}}}\right) \!\subseteq \sigma(\widehat A\,|\,\widehat \gI)$.
In view of Lemma \ref{lemma2}, $\gV_{\mathfrak p}$ is $(A,B,E,D_z)$-output nulling and therefore for all $x\in\gV_{\mathfrak p}$, we have $A_{F_{\mathfrak p}}\,x\in\gV_{\mathfrak p}$ and $(E+D_z\,F_{\mathfrak p})\,x=0$.
Using the closed-loop system equation we have $\widehat A\,\widehat{\gI}\subseteq\widehat{\gI}$, which implies $\mathfrak{p}(\widehat A\,\widehat{\gI})\subseteq\gV_{\mathfrak p}$. Thus, for each $\bsmat x\\[1mm] p\esmat\in \widehat \gI$
we have
$(A+B\,D_c\,W\,C)\,x+B\bigl(I+D_c(I-D_c\,D_y)^{-1}D_y\bigr)C_c\,p \in\gV_{\mathfrak p}$, 
which can be rewritten, by using the matrix inversion lemma, as
$(A+B\,D_c\,(I-D_y\,D_c)^{-1}\,C)\,x+(B\,\,(I-D_c\,D_y)^{-1}\,C_c)\,p \in\gV_{\mathfrak p}$. 
Combining the previous equation with $A_{F_{\mathfrak p}}\,x\in\gV_{\mathfrak p}$ yields
$B\big( (-F_{\mathfrak p}+\,D_c\,W\,C)x+(I-D_c\,D_y)^{-1}\,C_c\,p \big)\in\gV_{\mathfrak p}$. 
Moreover, since $\widehat{\gI} \subseteq \ker \widehat C$, we have 
$\widehat C \bsmat x \\[1mm] p \esmat=(E+D_z\,D_c\,W\,
C)\,x + (D_z\,C_c+D_z\,D_c\,W\,D_y\,C_c)\,p = 0$.
Combining the previous equation with $(E+D_z\,F_{\mathfrak p})\,x=0$ we obtain 
$D_z\big((-F_{\mathfrak p}+D_c\,W\,C)\,x + (C_c+D_z\,D_c\,W\,D_y\,C_c)\,p\big) = 0$.
Using again the matrix inversion lemma we find 
$D_z\big((-F_{\mathfrak p}+D_c\,W\,C)\,x + (I-D_c\,D_y)^{-1}C_c\,p\big) = 0$,
which implies $(-F_{\mathfrak p}+D_c\,W\,C)\,x + (I-D_c\,D_y)^{-1}C_c\,p\in\ker D_z$, so that
$B \big( (-F_{\mathfrak p}\!+\!D_c W C) x\!+\!(I\!-\!D_c D_y)^{-1} C_c p \big)\in B\ker D_z\cap\gV_{\mathfrak p}\subseteq\gR_{\gV_{\mathfrak p}}$. 
Adding and subtracting $A\,x$ to the left of the latter gives
$(A+B\,D_c\,W\,C)\,x+B\,(I-D_c\,D_y)^{-1}\,C_c\,\p-A_{F_{\mathfrak p}}\,\x \in \gR_{\gV_{\mathfrak p}}$,
form which
$\mathfrak{p}\,\left(\widehat{A}\,\bsmat \x \\[1mm] \p \esmat\right)-A_{F_{\mathfrak p}}\,\x\in \gR_{\gV_{\mathfrak p}}$. 
Let $\bsmat \x \\[1mm] \p \esmat$ be an eigenvector of $\widehat{A}$ in $\widehat{\gI}$ associated with the eigenvalue $\l$ and such that $\x\notin \gR_{\gV_{\mathfrak p}}$\footnote{Note that since we have $\gV_{\mathfrak p}=\mathfrak{p}\bigl(\widehat{\gI}\bigr)$, then $x\in \gV_{\mathfrak p}$. Thus, 
 if for every eigenvector $\bsmat \x \\[1mm] \p \esmat$ of $\widehat{A}|\widehat{\gI}$ we have $\x\in \gR_{\gV_{\mathfrak p}}$,
 then $\gR_{\gV_{\mathfrak p}}=\gV_{\mathfrak p}$, and in this case the claim is obvious because 
 $\sigma\left(A_{F_{\mathfrak p}}\,|\,\frac{\gV_{\mathfrak p}}{\gR_{\gV_{\mathfrak p}}}\right)=\varnothing$. 
 }, so that 
$\mathfrak{p}\,\left(\l\,\bsmat \x \\[1mm] \p \esmat\right)-A_{F_{\mathfrak p}}\,\x\in \gR_{\gV_{\mathfrak p}}$,
i.e., $\l\,\x-A_{F_{\mathfrak p}}\,\x\in \gR_{\gV_{\mathfrak p}}$ and $A_{F_{\mathfrak p}}\,\x\notin \gR_{\gV_{\mathfrak p}}$.
Consider the change of basis matrix $T=[\begin{array}{ccc} T_1 & T_2 & T_3 \end{array}]$ such that $T_1$ is a basis matrix of $\gR_{\gV_{\mathfrak p}}$, $[\begin{array}{ccc} T_1 & T_2  \end{array}]$ is a basis matrix of $\gV_{\mathfrak p}$. In this basis, partitioning $A_{F_{\mathfrak p}}$ and the vector $x$ conformably, the latter becomes
\[
\l \,\bmat{cc} \x_1 \\[-0.5mm]  \x_2 \\[-1.0mm]  0 \emat-\bmat{cccc}
A^{F_{\mathfrak p}}_{11} & A^{F_{\mathfrak p}}_{12} &A^{F_{\mathfrak p}}_{13} \\[-0.5mm] 
0 & A^{F_{\mathfrak p}}_{22} &A^{F_{\mathfrak p}}_{23} \\[-0.5mm] 
0 & 0 & A^{F_{\mathfrak p}}_{33} \emat\bmat{cc} \x_1 \\[-0.5mm]  \x_2 \\[-0.5mm]  0 \emat=\bmat{cc} \star \\[-0.5mm]  0 \\[-0.5mm]  0 \emat,
\]
from which we find $\l\,\x_2=A^{F_{\mathfrak p}}_{22}\,\x_2$. Since $\x_2 \neq 0$ (because $\x\notin \gR_{\gV_{\mathfrak p}}$), then $\l\in \sigma(A^{F_{\mathfrak p}}_{22})$. 
From the previous results, by defining $\Lambda$ as the set of eigenvalues of $\widehat A$ such that the corresponding eigenvector $\bsmat \x \\[1mm] \p \esmat$ satisfies $x\notin \gR_{\gV_{\mathfrak p}}$, we have 
$\Lambda\subseteq\sigma(A^{F_{\mathfrak p}}_{22})$. Since $\dim \widehat{\gI} \geq \dim \mathfrak{p}(\widehat{\gI})=\dim \gV_{\mathfrak p}$ we have ${\rm card}(\Lambda)\geq {\rm card}(\sigma(A^{F_{\mathfrak p}}_{22}))$, from which we obtain  $\sigma(A^{F_{\mathfrak p}}_{22})=\Lambda\subseteq \sigma(\widehat{A}\,|\,\widehat{\gI})\subseteq\sigma(\widehat{A})$. Finally,  by construction,  $\sigma(A^{F_{\mathfrak p}}_{22})=\sigma\left(A_{F_{\mathfrak p}}\,|\,\frac{\gV_{\mathfrak p}}{\gR_{\gV_{\mathfrak p}}}\right) \subseteq \sigma(\widehat A)$ for all  $F_{\mathfrak p}\in \mathfrak{F}_{\scriptscriptstyle (A,B,E,D_z)}(\gV_{\mathfrak p})$.
\endproof

The following result is a cornerstone of this paper: it shows that the set of matrices $K$ does not depend on the particular pair of subspaces $\gS$ and $\gV$ used to solve DDPDOF. 

\begin{theorem}
\label{lemK}
Let $\gS$ be $(A,H,C,G_y)$-self hidden and let $\gV$ be $(A,B,E,D_z)$-self bounded.
Let $(\gS,\gV;K)$ be a solution triple for DDPDOF. Then, for all $(A,H,C,G_y)$-self hidden subspaces $\bar{\gS}$ 
and all $(A,B,E,D_z)$-self bounded subspaces $\bar{\gV}$ that satisfy 
conditions {\bf (i-iii)} of Theorem~\ref{the66ths}, $(\bar{\gS},\bar{\gV};K)$ is a solution triple.
\end{theorem}
\proof
Consider the parameterization (\ref{param}) in Lemma~\ref{theparam}. We first prove that matrix $M$ is only dependent upon $\gS^\star_{\sy (A,H,C,G_y)}$: in other words, we want to show that $M$ can be chosen to be the same for any $(A,H,C,G_y)$-self hidden  
subspace $\gS$. Indeed, consider a basis matrix $T=[\begin{array}{ccc} T_1 & T_2 & T_3 \end{array}]$ of $\gS\oplus \gW$, such that $T_1$ is a basis matrix for $(\gS^\star_{\sy (A,H,C,G_y)}\oplus \gW)\cap \ker [\begin{array}{cc} C & G_y\end{array}]$ and $T_2$ extends $T_1$ to a basis for $\gS^\star_{\sy (A,H,C,G_y)}\oplus \gW$. From the $(A,H,C,G_y)$-self hiddenness of $\gS$, which can be expressed in terms of the inclusion  $\gS\oplus \gW \subseteq \gS^\star_{\sy (A,H,C,G_y)}\oplus \gW+\ker [\begin{array}{cc} C & G_y\end{array}]$, the columns of $T_3$ span a subspace of $\ker [\begin{array}{cc} C & G_y\end{array}]$. Thus, 
we can always choose $M=T_2$. Since both $T_1$ and $T_2$ can be chosen to be the same for any $(A,H,C,G_y)$-self hidden subspace $\gS$ (because they only depend on $\gS^\star_{\sy (A,H,C,G_y)}$ and $\ker [\begin{array}{cc} C & G_y\end{array}]$), matrix $M$ is only related to $\gS^\star_{\sy (A,H,C,G_y)}$. It follows, in particular, that if we choose $M$ to be the same for two different $(A,H,C,G_y)$-self hidden subspaces, then $[\begin{array}{cc} C & G_y\end{array}]\,M$ is the same, and such is also $\Psi$ in (\ref{param}).\\
We now show that $\Phi_2$ in (\ref{param}) does not depend on the particular 
$(A,B,E,D_z)$-self bounded subspace $\gV$ considered. Let $\bsmat \xi \\[1mm] \omega \esmat\in \ker \bsmat V && B \\[1mm]  0 && D_z\esmat$, where $V$ is a basis matrix for $\gV$. It suffices to show that there exists $\bar{\xi}$ such that $\bsmat \bar{\xi} \\[1mm] \omega \esmat\in \ker \bsmat \bar{V} && B \\[1mm]  0 && D_z\esmat$, where $\bar{V}$ is a basis matrix for $\bar{\gV}$. The condition 
$\bsmat \xi \\[1mm] \omega \esmat\in \ker \bsmat V && B \\[1mm]  0 && D_z\esmat$ can be written as 
$\bsmat V \\[1mm] 0 \esmat\,\xi=-\bsmat B \\[1mm] D_z \esmat\,\omega\in \gV\oplus 0_{\sy \gZ} \cap \ima \bsmat B \\[1mm] D_z \esmat$.
Since $\gV$ and $\bar{\gV}$ are $(A,B,E,D_z)$-self bounded, $\gV\oplus 0_{\sy \gZ} \cap \ima \bsmat B \\[1mm] D_z \esmat=\bar{\gV}\oplus 0_{\sy \gZ} \cap \ima \bsmat B \\[1mm] D_z \esmat$, so that $\bsmat V\\[1mm] 0 \esmat\,\xi\in \bar{\gV}\oplus 0_{\sy \gZ} \cap \ima \bsmat B \\[1mm] D_z \esmat$. Hence, there exists $\bar{\xi}$ such that $\bsmat  V\\[1mm] 0 \esmat\,\xi=\bsmat \bar{V} \\[1mm] 0 \esmat\,\bar{\xi}$. Thus, $\bsmat \bar{\xi} \\[1mm] \omega \esmat\in \ker \bsmat \bar{V} && B \\[1mm]  0 && D_z\esmat$.
Finally, we show that $R_2\,\tilde{A}\,M=R_2\,\bsmat A && H \\[1mm] E && G_z\esmat\,M$ in (\ref{param}) does not depend on the self bounded subspace and on the self hidden subspace considered. Let $\gS$ be $(A,H,C,G_y)$-self hidden satisfying $\ker [\begin{array}{cc} E & G_z \end{array}]\supseteq \left((\gS \oplus \gW) \cap \ker[\begin{array}{cc}C & G_y \end{array}]\right)$. From Theorem~\ref{prop3}, the smallest $(A,B,E,D_z)$-self bounded subspace $\gT$ that solves \\[-6mm]
\bea
\label{const1}
\ima \bmat{cc} H & S \\[-1.1mm] G_z & 0 \emat\subseteq \gT\oplus 0_{\sy \gZ} +\ima \bmat{cc} B \\[-1.1mm] D_z \emat
\eea
is $\gR^\star_{\sy (A,[\,B\;\;H\;\;S\,],E,[\,D_z\;\;G_z\;\;0])}$. The smallest $(A,B,E,D_z)$-self bounded subspace that solves DDPDOF together with the self hidden subspace $\gS$ is the smallest $(A,B,E,D_z)$-self bounded subspace $\gV^\prime$ such that 
\bea
\label{const2}
\ima  \bmat{cc} H \\[-1.1mm] G_z  \emat\subseteq \gV^\prime \oplus 0_{\sy \gZ} +\ima \bmat{cc} B \\[-1.1mm] D_z \emat \quad \text{and} \quad \gS\subseteq \gV^\prime.
\eea
Since (\ref{const2}) is a more stringent inclusion than (\ref{const1}), we have $\gV^\prime\supseteq \gT$. However, we show that $\gV^\prime$ and $\gT$ coincide. To this end, it suffices to show that $\gT$ satisfies (\ref{const2}). The fact that the first inclusion of (\ref{const2}) written for $\gT$ holds comes directly from (\ref{const1}). To prove that the second inclusion in (\ref{const2}) holds for $\gT$, i.e., that $\gS\subseteq \gT$, we use $\ima \bsmat S \\[1mm] 0 \esmat\subseteq \gT\oplus 0_{\sy \gZ}+\ima\bsmat B \\[1mm] D_z \esmat$, which is a consequence of 
(\ref{const1}), and rewrite it as $\gS\subseteq\gT +B\,\ker D_z$. We know also that $\gS \subseteq \gV^\prime$. Intersecting these inclusions gives
$\gS\subseteq \gV^\prime\cap (\gT+ B\,\ker D_z)=\gT+(\gV^\prime\cap B\,\ker D_z)=\gT$,
where we have used the modular rule (since $\gT\subseteq \gV^\prime$) and 
$\gV^\prime\cap B\,\ker D_z\subseteq \gV^\star_{\sy (A,B,E,D_z)} \cap B\,\ker D_z\subseteq \gT$.

Recall from Theorem~\ref{theparam} that the columns of $-R_2\,\tilde{A}\,M$ are the vectors $\u_i$ which, with suitable vectors $\xi_i$, satisfy
\bea
\label{eqfab}
\bmat{cc} A & H \\[-1.1mm] E & G_z \emat\,m_i=\bmat{cc} T\\[-1.1mm] 0 \emat\,\xi_i+\bmat{cc} B\\[-1.1mm] D_z \emat\,u_i
\eea
for all $i\in \{k+1,\ldots,r\}$, where $T$ is a basis matrix for $\gT$, $r$ is the dimension of $\gS \oplus \gW$ and $k$ is the dimension of $(\gS\oplus \gW)\cap \ker [\begin{array}{cc} C & G_y \end{array}]$, and $m_i$ is the $i$-th column of $M$. Solutions $\xi_i$ and $u_i$ to (\ref{eqfab}) exist because $\gT$ is $(A,B,E,D_z)$-output nulling and $\gT$ satisfies the solvability conditions of the decoupling problem. We now show that we can choose a different $(A,B,E,D_z)$-self bounded subspace containing $\gT$, say $\gV^\prime$ (with basis matrix $V^\prime$), and re-write (\ref{eqfab}) as 
\bea
\label{eqfab1}
\bmat{cc} A & H \\[-1mm] E & G_z \emat\,m_i=\bmat{cc} V^\prime \\[-1mm] 0 \emat\,\xi_i^\prime+\bmat{cc} B\\[-1mm] D_z \emat\,u_i
\eea
for some vectors $\xi_i^\prime$; notice that, as already shown, the vectors $m_i$ can be chosen to be the same for $\gT$ and $\gV^\prime$. Clearly, since $\gV^\prime$  contains $\gT$, the equation $T\,\xi_i=V^\prime\,\xi_i^\prime$ can always be solved in $\xi_i^\prime$.

We have shown that, given an $(A,H,C,G_y)$-self hidden subspace $\gS$, all the $(A,B,E,D_z)$-self bounded subspaces that, with $\gS$, form a solution to the decoupling problem have the same vectors $u_i$ (if, with no loss of generality, the $M$ matrices are chosen to be equal). Let us now consider in particular $\gS=\gS^\star_{\sy (A,H,C,G_y)}$, which is the infimum of all the $(A,H,C,G_y)$-self hidden subspaces, and we consider the subspace $\gT_{\rm min}=\gR^\star_{\sy (A,[\,B\;\;H\;\;S^\star\,],E,[\,D_z\;\;G_z\;\;0])}$, where $S^\star$ denotes a basis matrix for
$\gS^\star_{\sy (A,H,C,G_y)}$. The previous steps can be repeated verbatim for $\gT_{\rm min}$ to show that every other $(A,B,E,D_z)$-self bounded subspace which solves the problem has the same vectors $u_i$. Thus, $u_i$ (and therefore also the matrix $R_2\,\tilde{A}\,M$) are the same if $M$ is chosen to be the same. Finally, consider that (\ref{param}) parameterizes all possible matrices $K$ independently from the choice of $M$. This concludes the proof.
\endproof

We now show how to build a solution triple with a lower number of fixed poles.

\begin{lemma}
	\label{lemma2}
Let $(\gS,\gV;K)$  be a solution triple of DDPDOF,  and be such that $\gS$ is $(A,H,C,G_y)$-self hidden and $\gV$ is $(A,B,E,D_z)$-self bounded. Let 
$\bar{\gS}\defi \gS+(\gV_m\cap \gS_M)$.
Then, $(\bar{\gS},\bar{\gS}+\gV_m;K)$ is a solution triple satisfying 
\bea
\sigma_{\rm fixed}(\bar{\gS},\bar{\gS}+\gV_m)\subseteq \sigma_{\rm fixed}(\gS,\gV).
\eea
\end{lemma}
\proof 
Let $\bar{\gV}\defi \gV\cap (\gV_m+ \gS_M)$.
From Lemma~\ref{lem22}, $\gV_m+ \gS_M$ is $(A,B,E,D_z)$-self bounded and $\gV_m\cap \gS_M$ is $(A,H,C,G_y)$-self hidden. Thus, $\bar{\gV}$ is also $(A,B,E,D_z)$-self bounded and $\bar{\gS}$ is also $(A,H,C,G_y)$-self hidden.
We show that $(\bar{\gS},\bar{\gV};K)$ is a solution triple; since $(\gS,\gV;K)$ is a solution triple, 
$\ima \bsmat H \\[1mm] G_z \esmat\subseteq (\gV\oplus 0_{\scriptscriptstyle \gZ})+\ima \bsmat B \\[1mm] D_z\esmat$.
From $\ima \bsmat H \\[1mm] G_z \esmat\subseteq (\gV_m\oplus 0_{\scriptscriptstyle \gZ}) +\ima \bsmat B \\[1mm] D_z\esmat\subseteq (\gV_m+\gS_M) \oplus 0_{\scriptscriptstyle \gZ}+\ima \bsmat B \\[1mm] D_z\esmat$
and intersecting the previous two we obtain 
$\ima \bsmat H \\[1mm] G_z \esmat\subseteq (\bar{\gV} \oplus 0_{\scriptscriptstyle \gZ})+\ima \bsmat B \\[1mm] D_z\esmat$.
Using duality, $\ker [\begin{array}{cc} E & G_z \end{array}] \supseteq (\bar{\gS}\oplus \gW) \cap 
\ker [\begin{array}{cc} C & G_y \end{array}]$.
We show that $\bar{\gS}\subseteq \bar{\gV}$.

We prove that $\gV\supseteq \gV_m$. Since $\gV$ is $(A,B,E,D_z)$-self bounded and $\ima \bsmat H \\[1mm] G \esmat \subseteq (\gV \oplus 0_{\scriptscriptstyle \gZ}) +\ima \bsmat B \\[1mm] D_z\esmat$, then $\gV$ contains $\gV_m$, which is the smallest $(A,B,E,D_z)$-self bounded such that $\ima \bsmat H \\[1mm] G \esmat \subseteq (\gV_m \oplus 0_{\scriptscriptstyle \gZ}) +\ima \bsmat B \\[1mm] D_z\esmat$. Likewise, $\gS \subseteq \gS_M$. Since $\gS\subseteq \gV$ and $\gS \subseteq \gS_M\subseteq (\gS_M+\gV_m)$, then 
$\gS\subseteq \gV \cap(\gS_M+\gV_m) =\bar{\gV}$.
In a dual way, since $\gV \supseteq \gV_m$, then also $\gV\supseteq \gV_m\cap \gS_M$. Thus, also
$\gV_m\cap \gS_M\subseteq \gV\cap (\gV_m+ \gS_M)=\bar{\gV}$.
Summarizing, we found $\gS \subseteq \bar{\gV}$ and $\gV_m\cap \gS_M\subseteq	 \bar{\gV}$, which imply $\bar{\gS}=\gS+(\gV_m\cap \gS_M) \subseteq	 \bar{\gV}$.
Finally, Theorem~\ref{lemK} guarantees that $(\bar{\gS},\bar{\gV};K)$ is a solution triple.\\
We now show that $\sigma_{\rm fixed}(\bar{\gS},\bar{\gV})\subseteq \sigma_{\rm fixed}(\gS,\gV).$ 
We prove in particular that, given two solution triples $(\gS_1,\gV_1;K)$, 
 $(\gS_2,\gV_2;K)$ ($K$ can be the same from Theorem~\ref{lemK}), where $\gV_1$ and $\gV_2$ are $(A,B,E,D_z)$-self bounded subspaces with $\gV_1\subseteq \gV_2$ and such that $\gS_1$ and $\gS_2$ are $(A,H,C,G_y)$-self hidden subspaces with $\gS_1\supseteq \gS_2$, then $\sigma_{\rm fixed}(\gS_1,\gV_1)\subseteq \sigma_{\rm fixed}(\gS_2,\gV_2)$.

Since both $\gV$ and $\bar{\gV}$ contain $\gV_m$, and $\gV_m$ contains $\gR^\star_{\scriptscriptstyle (A,B,E,D_z)}$, 
\beann
\sigma\Big(A+B\,F_1\,\Big|\,\frac{\bar{\gV}\cap \langle A\,|\,\ima B \rangle}{\gR^\star_{\scriptscriptstyle (A,B,E,D_z)}}\Big)\subseteq 
\sigma\Big(A+B\,F_2\,\Big|\,\frac{{\gV}\cap \langle A\,|\,\ima B \rangle}{\gR^\star_{\scriptscriptstyle (A,B,E,D_z)}}\Big)
\eeann
for all $F_1\in \mathfrak{F}_{\scriptscriptstyle (A,B,E,D_z)}(\bar{\gV})$ and  $F_2\in \mathfrak{F}_{\scriptscriptstyle (A,B,E,D_z)}({\gV})$. 
Similarly
\beann
\sigma\Big(A\!+\!G_1\,C\,\Big|\,\frac{\gQ^\star_{\scriptscriptstyle (A,H,C,G_y)}}{\bar{\gS}\!+\!\langle \ker C\,|\,A \rangle}\Big)\!\subseteq 
\sigma\Big(A\!+\!G_2\,C\,\Big|\,\frac{\gQ^\star_{\scriptscriptstyle (A,H,C,G_y)}}{\gS\!+\!\langle \ker C\,|\,A \rangle}\Big)
\eeann
for all $G_1\in \mathfrak{G}_{\scriptscriptstyle (A,H,C,G_y)}(\bar{\gS})$ and $G_2\in \mathfrak{G}_{\scriptscriptstyle (A,H,C,G_y)}({\gS})$. Since both $\gV$ and $\bar \gV$ are $(A,B,E,D_z)$-self bounded and $\gV \supseteq \bar \gV$, we have $\sigma_{\rm fixed}(\bar{\gV})\subseteq \sigma_{\rm fixed}(\gV)$\footnote{This follows from the definition of fixed poles and from the fact that self bounded subspaces share the same reachability subspace}. Similarly, since 
both $\gS$ and $\bar \gS$ are $(A,H,C,G_y)$-self hidden and $\gS \subseteq \bar \gS$, we find that $\sigma_{\rm fixed}(\bar{\gS})\subseteq
\sigma_{\rm fixed}(\gS)$, which implies
$\sigma_{\rm fixed}(\bar{\gS},\bar{\gV})\subseteq \sigma_{\rm fixed}({\gS},{\gV})$.
Now we show that $(\bar{\gS},\bar{\gS}+\gV_m;K)$ is a solution triple with a smaller number of fixed poles. From the solvability of DDPDOF, we have
 \bea
 \label{pk1}
 \ima \bmat{cc} H \\[-1mm] G_z \emat \subseteq  (\gV^\star_{\scriptscriptstyle (A,B,E,D_z)} \oplus {0}_{\gZ}) +\ima \bmat{cc} B \\[-1mm] D_z\emat,
 \eea
 and, moreover, $\bar{\gS}\subseteq \bar{\gV}\subseteq \gV^\star_{\scriptscriptstyle (A,B,E,D_z)}$. We can write these two inclusions together as
 \bea
 \label{pk2}
  \ima \bmat{cc} H & \bar{S} \\[-1mm] G_z & 0 \emat \subseteq  (\gV^\star_{\scriptscriptstyle (A,B,E,D_z)} \oplus {0}_{\gZ}) +\ima \bmat{cc} B \\[-1mm] D_z\emat,
 \eea
where $\bar{S}$ is a basis matrix of $\bar{\gS}$. We can now apply 
Theorem~\ref{prop3}, which shows that 
$\gR^\star_{\scriptscriptstyle (A,[\,B\;\;H\;\;\bar{S}\,],E,[\,D_z\;\;G_z\;\;0\,])}$ is the smallest 
$(A,B,E,D_z)$-self bounded subspace  that satisfies 
$\ima \bsmat H & \bar{S} \\[1mm] G_z & 0 \esmat \subseteq   (\gR^\star_{\scriptscriptstyle (A,[\,B\;\;H\;\;\bar{S}\,],E,[\,D_z\;\;G_z\;\;0\,])} \oplus {0}_{\gZ}) +\ima \bsmat B \\[1mm] D_z\esmat$.
Eq. (\ref{pk1}) holds also with $\bar{\gV}$ in place of $\gV^\star_{\scriptscriptstyle (A,B,E,D_z)}$, and also $\bar{\gS}\subseteq \bar{\gV}$; thus $\bar{\gV}$ also satisfies (\ref{pk2}) with $\bar{\gV}$ in place of $\gV^\star_{\scriptscriptstyle (A,B,E,D_z)}$. The infimality of $\gR^\star_{\scriptscriptstyle (A,[\,B\;\;H\;\;\bar{S}\,],E,[\,D_z\;\;G_z\;\;0\,])}$ ensures that $\gR^\star_{\scriptscriptstyle (A,[\,B\;\;H\;\;\bar{S}\,],E,[\,D_z\;\;G_z\;\;0\,])}\subseteq \bar{\gV}$, so that $\sigma_{\rm fixed}(\bar{\gS},\gR^\star_{\scriptscriptstyle (A,[\,B\;\;H\;\;\bar{S}\,],E,[\,D_z\;\;G_z\;\;0\,])})\subseteq \sigma_{\rm fixed}(\bar{\gS},\bar{\gV})$. 
Using Lemma~\ref{lemA2m}, we have $\gR^\star_{\scriptscriptstyle (A,[\,B\;\;H\;\;\bar{S}\,],E,[\,D_z\;\;G_z\;\;0\,])}=\bar{\gS}+\gV_m$, and this concludes the proof.
\endproof

\begin{lemma}
	\label{lem_init}
Let DDPDOF be solvable. 
Let the pair $(A,[\begin{array}{cc} B & H \end{array}])$ be reachable and let the pair $\left(\bsmat C\\[1mm] E \esmat,A\right)$ be observable.
Let $\sigma^\ast\defi \sigma_{\rm fixed}(\gV_m\cap \gS_M,\gV_m)$. Then, there exists $K$ such that for every $(A,H,C,G_y)$-input containing subspace $\gS$ such that $\gV_m\cap \gS_M\subseteq \gS\subseteq \gS_M$, the triple $(\gS,\gS+\gV_m;K)$ is a solution triple and
$\sigma_{\rm fixed}(\gS,\gS+\gV_m)=\sigma^\ast$.
\end{lemma}
\proof
From Property~\ref{pA4}, the subspace $\gS+\gV_m$ is $(A,B,E,D_z)$-self bounded, which implies that $\gS+\gV_m$ is $(A,B,E,D_z)$-output nulling. By definition, $\gS$ is $(A,H,C,G_y)$-input containing, and it is contained in $\gS+\gV_m$. 
 Finally, Theorem \ref{the66ths} and Theorem~\ref{lemK} ensure the existence of $K$ such that $(\gS,\gS+\gV_m;K)$ is a solution triple. Its fixed spectrum is
$\sigma_{\rm fixed}(\gS,\gS+\gV_m)=\sigma_{\rm fixed}(\gS)\uplus \sigma_{\rm fixed}(\gS+\gV_m)$.
From the reachability of the pair $(A,[\begin{array}{cc} B & H \end{array}])$, we obtain $\langle A\,|\,\ima B +\ima H \rangle=\gX$. In view of Lemma~\ref{A5m}, this implies that $\langle A\,|\,\ima B \rangle+\gV_m=\gX$, which in turn implies
$\langle A\,|\,\ima B \rangle+\gS+\gV_m=\gX$.
 Dually, 
the observability of the pair $\left(\bsmat C\\[1mm] E \esmat,A\right)$ implies $\langle \ker C \cap \ker E\,|\,A\rangle=0_{\scriptscriptstyle \gX}$, so that, by Lemma~\ref{A5m}, we have $\langle \ker C\,|\,A\rangle\cap \gS_M=0_{\scriptscriptstyle \gX}$, which in turn implies
$\langle \ker C\,|\,A\rangle\cap \gS=0_{\scriptscriptstyle \gX}$,
since $\gS\subseteq \gS_M$.
Let us now consider $\sigma_{\rm fixed}(\gS+\gV_m)$. Defining $A_F=A+B\,F$ and $A_G=A+G\,C$, we have by definition of fixed poles
\bea
&&\sigma_{\rm fixed}(\gS+\gV_m) \nonumber  \\
\ns&\ns \ns&\ns \hspace{4mm} =
\sigma\Big(A_F\,\Big|\,\frac{\gS+\gV_m}{\gR_{\gS+\gV_m}}\Big)
\uplus \sigma\Big(A_F\,\Big|\,\frac{\gX}{\gS+\gV_m+\langle A\,|\,\ima B \rangle}\Big) \nonumber \\
\ns&\ns \ns&\ns \hspace{4mm} =
\sigma\Big(A_F\,\Big|\,\frac{\gS+\gV_m}{\gR^\star_{\scriptscriptstyle (A,B,E,D_z)}}\Big)
\uplus \sigma\Big(A_F\,\Big|\,\frac{\gX}{\gX}\Big)=
\sigma\Big(A_F\,\Big|\,\frac{\gS+\gV_m}{\gR^\star_{\scriptscriptstyle (A,B,E,D_z)}}\Big)  \nonumber \\
\ns&\ns \ns&\ns \hspace{4mm} =
\sigma\Big(A_F\,\Big|\,\frac{\gS+\gV_m}{\gV_m}\Big)\uplus
\sigma\Big(A_F\,\Big|\,\frac{\gV_m}{\gR^\star_{\scriptscriptstyle (A,B,E,D_z)}}\Big) \nonumber  \\
\ns&\ns \ns&\ns \hspace{4mm} =
\sigma\Big(A_G\,\Big|\,\frac{\gS}{\gV_m\cap \gS_M}\Big)\uplus
\sigma\Big(A_F\,\Big|\,\frac{\gV_m}{\gR^\star_{\scriptscriptstyle (A,B,E,D_z)}}\Big), \nonumber 
\eea
where we have used 
Lemma~\ref{A3M}, and 
the fact that $\gR_{\gS+\gV_m}=\gR^\star_{\scriptscriptstyle (A,B,E,D_z)}$. Indeed, since $\gS+\gV_m$ is $(A,B,E,D_z)$-self bounded, it contains $\gR^\star_{\scriptscriptstyle (A,B,E,D_z)}$.
Similarly, it is seen that
\bea
 \sigma_{\rm fixed}(\gS) 
=\sigma\Big(A_G\,\Big|\,\frac{\gQ^\star_{\scriptscriptstyle (A,H,C,G_y)}}{\gS_M}\Big)\uplus
\sigma\Big(A_G\,\Big|\,\frac{\gS_M}{\gS}\Big).
\eea
Therefore,
\beann
&&\sigma_{\rm fixed}(\gS,\gS+\gV_m) \\
\ns&\ns \ns&\ns \hspace{4mm} =
 \sigma_{\rm fixed}(\gS)\uplus
\sigma_{\rm fixed}(\gS+\gV_m) =
 \sigma\Big(A_G\,\Big|\,\frac{\gQ^\star_{\scriptscriptstyle (A,H,C,G_y)}}{\gS_M}\Big) \\
\ns&\ns \ns&\ns \hspace{4mm} \uplus
\sigma\Big(A_G\,\Big|\,\frac{\gS_M}{\gS}\Big) \uplus
\sigma\Big(A_G\,\Big|\,\frac{\gS}{\gV_m\cap \gS_M}\Big)  \uplus\,
\sigma\Big(A_F\,\Big|\,\frac{\gV_m}{\gR^\star_{\scriptscriptstyle (A,B,E,D_z)}}\Big) \\
\ns&\ns \ns&\ns \hspace{4mm} =
 \sigma\Big(\!A_G\,\Big|\,\frac{\gQ^\star_{\scriptscriptstyle (A,H,C,G_y)}}{\gV_m\cap \gS_M}\Big)\!\uplus
\sigma\Big(\!A_F\,\Big|\,\frac{\gV_m}{\gR^\star_{\scriptscriptstyle (A,B,E,D_z)}}\Big)\\
\ns&\ns \ns&\ns \hspace{4mm} =
\sigma_{\rm fixed}(\gV_m\cap \gS_M)\uplus \sigma_{\rm fixed}(\gV_m).\\[-1.2cm]
\eeann
\endproof

The previous lemma (see also \cite[Lemma 4]{DMCuellar-M-01}) says that for every input containing between $\gV_m\cap \gS_M$ and $\gS_M$, the triple $(\gS,\gS+\gV_m;K)$ is a solution triple with the same fixed poles. This gives rise to a family of  solution triples with the same fixed poles $\sigma^\ast$:
\bea
&&\sigma_{\rm fixed}(\gV_m\cap \gS_M,\gV_m) \qquad \text{(taking $\gS=\gV_m\cap \gS_M$)}, \nonumber \\[-1.0mm]
&&\qquad \vdots \label{chain1} \\[-1.0mm]
&& \sigma_{\rm fixed}(\gS_M, \gS_M+\gV_m)\qquad \text{(taking $\gS=\gS_M$).} \nonumber 
\eea
The following is the dual of the latter.

\begin{lemma}
Let DDPDOF be solvable. Let the pair $(A,[\begin{array}{cc} B & H \end{array}])$ be reachable and let the pair $\left(\bsmat C\\[1mm] E \esmat,A\right)$ be observable.
Let $\sigma^\ast\defi \sigma_{\rm fixed}(\gS_M,\gS_M+\gV_m)= \sigma_{\rm fixed}(\gV_m\cap\gS_M,\gV_m)$. Then, for every $(A,B,E,D_z)$-output nulling subspace $\gV$ such that $\gV_m\subseteq \gV\subseteq \gV_m+\gS_M$, there exists $K$ such that $(\gV\cap \gS_M,\gV;K)$ is a solution triple and 
$\sigma_{\rm fixed}(\gV\cap \gS_M,\gV)=\sigma^\ast$.
\end{lemma}

The previous lemma says that for every output nulling subspace between $\gV_m$ and $\gS_M+\gV_m$, the triple $(\gV\cap \gS_M,\gV; K)$ is a solution triple with the same fixed poles. This gives rise to a family of solution triples with the same fixed poles $\sigma^\ast$:
\bea
&&\sigma_{\rm fixed}(\gV_m\cap \gS_M,\gV_m) \qquad \text{(taking $\gV=\gV_m$),}\nonumber  \\[-1.0mm]
&& \qquad \vdots \label{chain2}\\[-1.0mm]
&& \sigma_{\rm fixed}(\gS_M, \gS_M+\gV_m)\qquad \text{(taking $\gV=\gV_m+\gS_M$).} \nonumber 
\eea
Clearly, this set of solutions is exactly the same that was obtained before.
We now eliminate the assumption of observability of the pair $\left(\bsmat C\\[1mm] E \esmat,A\right)$.
\begin{corollary}
Let DDPDOF be solvable. 
Let the pair $(A,[\begin{array}{cc} B & H \end{array}])$ be reachable.
Let $\sigma^\ast\defi \sigma_{\rm fixed}(\gV_m\cap \gS_M,\gV_m)$. Then, for every $(A,H,C,G_y)$-input containing subspace $\gS$ such that $\gV_m\cap \gS_M\subseteq \gS\subseteq \gS_M$, there exists $K$ such that $(\gS,\gS+\gV_m;K)$ is a solution triple and
$\sigma_{\rm fixed}(\gS,\gS+\gV_m)\supseteq \sigma^\ast$.
\end{corollary}
\proof
Since $\gV_m\cap \gS_M\cap \langle \ker C\,|\,A\rangle\subseteq \gS\cap \langle \ker C\,|\,A\rangle$,
clearly 
\[
\sigma^\ast=\sigma_{\rm fixed}(\gV_m\cap \gS_M,\gV_m) \subseteq \sigma_{\rm  fixed}(\gS,\gS+\gV_m).
\]
\endproof

When we drop the observability assumption, it is no longer true that all the pairs in (\ref{chain1})
 have the same number of fixed poles: the ``smallest'' one, i.e., 
$\sigma_{\rm fixed}(\gV_m\cap \gS_M,\gV_m)$, minimizes the number of fixed poles, and the ``larger'' the pair in the list, the greater the number of elements in $\sigma_{\rm fixed}$.
We now drop the reachability assumption.
\begin{corollary}
	\label{cor_fin}
Let DDPDOF be solvable. 
Let the pair $\left(\bsmat C\\[1mm] E \esmat,A\right)$ be observable.
Let $\sigma^\ast \defi \sigma_{\rm fixed}(\gS_M,\gS_M+\gV_m)$. For every $(A,B,E,D_z)$-output nulling subspace $\gV$ such that $\gV_m\subseteq \gV\subseteq \gV_m+\gS_M$, there exists $K$ such that $(\gV\cap \gS_M,\gV;K)$ is a solution triple and
$\sigma_{\rm fixed}(\gV\cap \gS_M,\gV)\supseteq\sigma^\ast$.
\end{corollary}

\proof 
Since $\gV_m+\gS_M \supseteq \gV$, it follows that $\sigma^\ast=\sigma_{\rm fixed}(\gS_M,\gS_M+\gV_m)\subseteq \sigma_{\rm  fixed}(\gV\cap \gS_M,\gV)$.
\endproof

When we drop the reachability assumption, it is no longer true that all the pairs in (\ref{chain2})
have the same number of fixed poles: the ``largest'' one, i.e., 
$\sigma_{\rm fixed}(\gS_M, \gS_M+\gV_m)$, minimizes the number of fixed poles, and the ``smallest'' the pair in the list, the greater the number of $\sigma_{\rm fixed}$.

\begin{remark}
\label{rem61}
	{
		A key contribution of this paper is to show that the well-posedness condition is decoupled from the problem of the fixed poles. This is a direct consequence of Theorem~\ref{lemK}. In fact, if DDPDOF is solvable, the set $\gK$ of matrices such $(\gS,\gV;K)$ is a solution triple, where $\gS$ is $(A,H,C,G_y)$-self hidden and $\gV$ is $(A,B,E,D_z)$-self bounded, coincides with set $\gK^\star$ of matrices such that
		\beann
		\bmat{cc} A+B\,K\,C & H+B\,K\,G_y \\
		E+D_z\,K\,C & G_z+ D_z\,K\,G_y \emat \!\!(\gS^\star_{\sy (A,H,C,G_y)}  \oplus \gW)\subseteq \gV^\star_{\sy (A,B,E,D_z)} \oplus 0_{\scriptscriptstyle \scriptscriptstyle \gZ}.  
		\eeann
		Indeed, from the minimality of $\gS^\star_{\sy (A,H,C,G_y)}$ and the maximality of $\gV^\star_{\sy (A,B,E,D_z)}$ we have 
		$\gK \subseteq \gK^\star$ (see \cite{PN-submitted}). On the other hand, from the solvability of the problem, for all $K \in \gK^\star$, $(\gS^\star_{\sy (A,H,C,G_y)},\gV^\star_{\sy (A,B,E,D_z)};K)$ is a solution triple. Since  $\gS^\star_{\sy (A,H,C,G_y)}$ is $(A,H,C,G_y)$-self hidden and $\gV^\star_{\sy (A,B,E,D_z)}$ is $(A,B,E,D_z)$-self bounded, Theorem~\ref{lemK} ensures that $(\gS,\gV;K)$ is also a solution triple, so that $\gK^\star\subseteq\gK$ which gives $\gK=\gK^\star$.  
		}
\end{remark}

\begin{theorem}
	\label{the6.8}
	Let DDPDOF be solvable and let either the pair $(A,[\begin{array}{cc} B & H \end{array}])$ be reachable or the pair $\left(\bsmat C\\[1mm] E \esmat,A\right)$ be observable or both. 
Let $\aleph \defi \{\sigma_{\rm fixed}(\gS,\gV)\,|\,\exists K:$ $(\gS,\gV;K)$ is a solution triple with $\gS$ $(A,H,C,G_y)$-self hidden and $\gV$ $(A,B,E,D_z)$-self bounded$\}$. Then, $\aleph$ 
	has a minimal element $\sigma^\star$ satisfying
	\bea
	\label{sigmastar}
	\sigma^\star\!=\min \{\sigma_{\rm fixed}(\gS_M,\gS_M\!+\!\gV_m),
		\;\sigma_{\rm fixed}(\gV_m\cap \gS_M,\gV_m) \}
	\eea
	 and $\sigma^\star\subseteq \sigma(\widehat{A})$ for every controller that solves DDPDOF.
\end{theorem}

\proof
First, from Lemma~\ref{giapponesi}, if DDPDOF is solvable, a solution triple $(\gS,\gV;K)$ exists where $\gS$ is $(A,H,C,G_y)$-self hidden  and $\gV$ is $(A,B,E,D_z)$-self bounded such that $\sigma_{\rm  fixed}(\gS,\gV)\subseteq\sigma(\widehat{A})$. Thus, $\aleph$ is non-empty, and if $\aleph$ admits minimum, the last claim is proved.
Now, we show that $\aleph$ has a minimal element, and it is exactly $\sigma^\star$. Since either $\sigma_{\rm fixed}(\gS_M,\gS_M+\gV_m) \supseteq
\sigma_{\rm fixed}(\gV_m\cap \gS_M,\gV_m) $ if the pair $(A,[\begin{array}{cc} B & H \end{array}])$ is reachable, or $\sigma_{\rm fixed}(\gV_m\cap \gS_M,\gV_m) \supseteq \sigma_{\rm fixed}(\gS_M,\gS_M+\gV_m)$ if the pair $\left(\bsmat C\\[1mm] E \esmat,A\right)$ is  observable, then  $\min(\cdot)$ operation in \eqref{sigmastar} is well-defined.
In particular, assume that a solution triple $(\gS,\gV;K)$ exists where $\gS$ is $(A,H,C,G_y)$-self hidden and $\gV$ is $(A,B,E,D_z)$-self bounded such that $\hat{\sigma}=\sigma_{\rm  fixed}(\gS,\gV)\subseteq\sigma^\star$. We show that $\sigma^\star=\hat{\sigma}$. Define $\bar{\gS}=\gS+(\gV_m\cap \gS_M)$. From Lemma \ref{lemma2} and Theorem~\ref{lemK}, $(\bar{\gS},\,\bar{\gS}+\gV_m;K)$ is a solution triple and 
$\sigma_{\rm fixed}(\bar{\gS},\,\bar{\gS}+\gV_m)\subseteq \sigma_{\rm fixed}(\gS,\gV)=\hat{\sigma}$. Since $(\gV_m\cap \gS_M) \subseteq \bar{\gS} \subseteq \gS_M$  (where the last inclusion comes from Theorem~\ref{the12d3}), in view of Corollary \ref{cor_fin} - Lemma \ref{lem_init}  we have $\sigma_{\rm fixed}(\bar{\gS},\,\bar{\gS}+\gV_m)\supseteq \sigma^\star$, which implies $\sigma^\star \subseteq \hat{\sigma}$. The latter, together with $\hat{\sigma} \subseteq \sigma^\star $ gives  $\hat{\sigma} = \sigma^\star$. 
\endproof

Finally, we address the case where we allow $(A,[\begin{array}{cc} B & H \end{array}])$ to be non-reachable and $\left(\bsmat C\\[1mm] E \esmat,A\right)$ to be non-observable.
In this case, we can only establish upper and lower bounds for the set of fixed poles, which might no longer admit a minimum. 

\begin{theorem}
		Let DDPDOF be solvable. Define
	\beann
	\sigma^\ddagger \defi \sigma_{\rm fixed}(\gS_M\!\cap\! \gV_m,\gS_M\!+\!\gV_m) \quad\! \!\!\!\text{and}\!\!\! \!\quad
	\sigma^\dagger \defi \sigma_{\rm fixed}(\gS_M,\gV_m).
	\eeann
	The set $\aleph$
	has elements $(\bar \gS,\bar \gV)$ satisfying
	$\sigma_{\rm fixed}(\bar{\gS},\bar{\gV})\subseteq\sigma^\ddagger$
	and $\sigma^\dagger\subseteq \sigma(\widehat{A})$ for every controller solving DDPDOF.
\end{theorem}
\proof Let $(\gS,\gV;K)$ be a solution triple for DDPDOF where $\gS$ is $(A,H,C,G_y)$-self hidden and $\gV$ is $(A,B,E,D_z)$-self bounded: this is not restrictive in view of Lemma~\ref{giapponesi}.
As in the proof of Theorem~\ref{the6.8}, we define $\bar{\gS}=\gS+(\gV_m\cap \gS_M)$, and observe that 
 $(\bar{\gS},\,\bar{\gS}+\gV_m;K)$ is a solution triple and 
$\sigma_{\rm fixed}(\bar{\gS},\,\bar{\gS}+\gV_m)\subseteq \sigma_{\rm fixed}(\gS,\gV)$. By definition $\sigma_{\rm fixed}(\bar{\gS},\bar{\gS}+\gV_m)=\sigma_{\rm fixed}(\bar{\gS})\uplus\sigma_{\rm fixed}(\bar{\gS}+\gV_m)$, and the first claim follows on defining $\bar \gV=\bar{\gS}+\gV_m$ by noting that $\sigma_{\rm fixed}(\bar{\gS})\subseteq \sigma_{\rm fixed}({\gS_M\cap \gV_m})$ and $\sigma_{\rm fixed}(\bar \gV)=\sigma_{\rm fixed}(\bar{\gS}+\gV_m)\subseteq \sigma_{\rm fixed}({\gS_M+ \gV_m})$.
The second claim is proved by contradiction. Assume that there exists a controller that solves DDPDOF and such that $\sigma^\dagger\cap \sigma(\widehat{A})\neq \sigma^\dagger$. In view of Lemma~\ref{giapponesi} there exists a solution triple $(\gS,\gV;K)$  of DDPDOF where $\gS$ is $(A,H,C,G_y)$-self hidden  and $\gV$ is $(A,B,E,D_z)$-self bounded such that $\sigma_{\rm  fixed}(\gS,\gV)\subseteq\sigma(\widehat{A})$.
Following the first part of this proof, we find that $(\bar{\gS},\,\bar{\gV};K)$ is also a solution triple and 
$\sigma_{\rm fixed}(\bar{\gS},\,\bar{\gV})\subseteq \sigma_{\rm fixed}(\gS,\gV)\subseteq\sigma(\widehat{A})$. However, $\sigma_{\rm fixed}(\bar \gS)\supseteq \sigma_{\rm fixed}({\gS_M})$ and $\sigma_{\rm fixed}(\bar{\gV})\supseteq \sigma_{\rm fixed}({\gV_m})$ so that $\sigma_{\rm fixed}(\bar{\gS},\,\bar{\gV})\supseteq \sigma^\dagger$, which implies that $\sigma^\dagger=\sigma^\dagger\cap \sigma_{\rm  fixed}(\bar \gS,\bar \gV)
\subseteq\sigma^\dagger \cap \sigma(\widehat{A})$, i.e, $\sigma^\dagger=\sigma^\dagger \cap \sigma(\widehat{A})$ leading to contradiction.
\endproof

\begin{example}
{
Consider the following 
system:
\beann
A\ns&\ns =\ns&\ns 
\bsmat
-30 && 0 && 0 && 0 \\[1mm]
0 && 0 && 0 && 0 \\[1mm]
0 && 0 && -1 && 0 \\[1mm]
0 && 0 && 0 && 0
\esmat, \;\; B=
\bsmat
0 && 1 && 0 \\[1mm]
10 && 0 && 0 \\[1mm]
0 && 13 && -1 \\[1mm]
0 && 0 && 1
\esmat,\;\; H=\bsmat
-1 \\[1mm]
0 \\[1mm]
0 \\[1mm]
-1
\esmat,\\
C\ns&\ns =\ns&\ns 
\bsmat
0 && -1 && 0 && -1 \\[1mm]
0 && 0 && 0 && 0 \\[1mm]
-14 && 0 && 0 && 1
\esmat,\;\; D_y=
\bsmat
-11 && 0 && 0 \\[1mm]
-1 && 13 && -5 \\[1mm]
0 && 0 && -1
\esmat, \;\; G_y=
\bsmat
-5 \\[1mm]
-1 \\[1mm]
0
\esmat, \\
E\ns&\ns =\ns&\ns\bsmat
0 && 0 && 0 && -20
\esmat,\;\; D_z=
\bsmat
0 && 0 && -1
\esmat,\;\; G_z=1.
\eeann
It is easy to see that $Z_{(A,B,E,D_z)}=\{-20\}$ and $Z_{(A,H,C,G_y)}=\{-1\}$, and $\gV^\star_{\sy (A,B,E,D_z)}=\gX$ and $\gS^\star_{\sy (A,H,C,G_y)}=\{0\}$. The conditions of Lemma 5.2 in \cite{PN-submitted} guarantee the existence of a matrix $K$ such that (\ref{clK}) holds. 
Using Theorem~\ref{theparam}, the set of all matrices $K$ that solve the problem is parameterized as
$K=\bsmat \star && \star &&\star \\[0mm]
 \star && \star &&\star \\[0mm]
k_{31} && k_{32} &&\star \esmat$, 
where $\star$ denote arbitrary values and $5\,k_{31}+k_{32}=-1$. Choosing for simplicity all the arbitrary values to be zero, the only value of $k_{31}$ such that $I+K\,D_y$ is singular is $k_{31}=-6/25$. We choose for example $k_{31}=1$ and $k_{32}=-6$. It is easy to see that the triple $(\gS^\star_{\sy (A,H,C,G_y)},\gV^\star_{\sy (A,B,E,D_z)};K)$ solves the problem.
The compensator built from $\gV^\star_{\sy (A,B,E,D_z)}$ and $\gS^\star_{\sy (A,H,C,G_y)}$ leads to a closed-loop spectrum in the form $\{-20,-1\}\cup \sigma_{\rm free}$, where $\sigma_{\rm free}$ is completely assignable with a suitable choice of the matrices $F$ and $G$. For example
\[
F=\bsmat 0 && -0.35 && 0 && 0 \\[1mm]
27.03 && 0 && 0 && -1.52 \\[1mm]
0 && 0&& 0&& -20\esmat
\;  \;\text{and}\;\;
G=\bsmat
-0.93 && 3.67 && -1.87 \\[1mm]
1.81 && -9.06 && 0.13 \\[1mm]
0 && 0 && 0 \\[1mm]
1.69 && -9.44 && -0.13\esmat
   \]
   lead to the closed-loop spectrum $\{-1,-20,-1.5,-2.5,-3.5\}$, where the last 3 have double multiplicity. In this case
   \beann
   \gV_m=\ima \bsmat 1 && 0 && 0 \\ 0 && 1 && 0 \\ 0 && 0 && 1 \\ 0 && 0 && 0 \esmat\quad \text{and} \quad \gV_m\cap \gS_M=\bsmat 0\\0\\1\\0 \esmat.
   \eeann
   Notice that the subspaces $\gV^\star_{\sy (A,B,E,D_z)}$, $\gS^\star_{\sy (A,H,C,G_y)}$, $\gV_m$ and $\gS_M$ can all be obtained with the standard sequences (6) and (7) in \cite{PN-submitted}.
   Theorem~\ref{lemK} ensures that the same matrices $K$ solve the problem with $\gV_m$ and $\gV_m\cap \gS_M$; thus we take the same $K$ as in the previous case (the closed-loop spectrum does not depend on $K$).
Solving the problem using the compensator built from these subspaces leads to a closed-loop spectrum in the form $\{-1\}\cup \sigma_{\rm free}$, where $\sigma_{\rm free}$ is completely assignable with a suitable choice of the matrices $F$ and $G$. For example
$F=\bsmat 0 && -0.35 && 0 && 0 \\[1mm]
27.03 && 0 && 0 && 0 \\[1mm]
0 && 0 && 0 && -1.50\esmat$ and the same matrix $G$ as above
   lead to the closed-loop spectrum $\{-1,-1.5,-2.5,-3.5\}$, where $-1.5$ is triple and $-2.5$ and $-3.5$ are double, thus eliminating the high frequency zero from the  spectrum. 
   }
\end{example}

\begin{example}
{
Consider the system described by
\beann
A\ns&\ns =\ns&\ns 
\bsmat
-1 && 0 && 0 && 0 \\[1mm]
7 && 0 && -6 && 0 \\[1mm]
0 && -2 && 0 && 0 \\[1mm]
0 && 0 && -9 && 0
\esmat, \;\; B=\bsmat
-9 && 0 && 0 \\[1mm]
-1 && 1 && 10 \\[1mm]
0 && 0 && 0 \\[1mm]
0 && -6 && 0
\esmat, \;\; H=\bsmat
0 && 0 \\[1mm]
-8 && 0 \\[1mm]
0 && 7 \\[1mm]
0 && 3
\esmat, \\
C\ns&\ns =\ns&\ns 
\bsmat
0 && 8 && 1 && -10 \\[1mm]
0 && 0 && 8 && -1
\esmat,\;\; 
D_y=\bsmat
-5 && 0 && 0 \\[1mm]
0 && -1 && -2
\esmat,\;\; G_y=\bsmat
0 && 2 \\[1mm]
0 && 0
\esmat,\\
E\ns&\ns =\ns&\ns \bsmat
29 && 0 && 0 && 0
\esmat, \;\; D_z=
\bsmat
-9 && 0 && 0
\esmat,\;\; G_z=
\bsmat
0 && 0
\esmat.
\eeann
Here $Z_1=\{-30\}$ and $Z_2=\{-1,-11/76\}$.
In this case $\gS^\star_{\sy (A,H,C,G_y)} =\spanR\left\{\bsmat 0\\1\\0\\0\esmat,\bsmat 0\\0\\5\\2 \esmat\right\}$ and $\gV^\star_{\sy (A,B,E,D_z)}=\gX$. In this case the problem is solvable with $K=0$. However, using $\gV^\star_{\sy (A,B,E,D_z)}$ and $\gS^\star_{\sy (A,H,C,G_y)}$ forces the closed-loop spectrum to contain all the zeros $Z_1\cup Z_2$, and in particular a high frequency mode $-30$ and a low frequency mode at $-11/76$. We solve the problem using e.g. the triple $(\gS_M,\gV_m+\gS_M;K)$. In this case $\gV_m+\gS_M=\gS_M=\ima \bsmat 0 && 0 && 0 \\ 1 && 0 && 0\\ 0 && 1 && 0 \\ 0 && 0 && 1 \esmat$. Solving the problem with these subspaces and $K=0$ is such that the only fixed pole in the closed-loop spectrum is $-1$. In other words, this solution eliminates $-30$ and $-11/76$ from the closed-loop spectrum.  The same result is obtained taking e.g. $(\gV_m\cap \gS_M,\gV_m;K)$.
}
\end{example}

\section*{Appendix A}
We recall some fundamental results on geometric control theory, most of which are restatements or dualizations of the results in \cite[Appx.~A]{DMCuellar-M-01} and \cite[Lemma 3]{N-TAC-08}.
We first consider the inclusion $\ima L \subseteq \gV^\star_{\scriptscriptstyle (A,B,C,D)}$, which is the solvability condition of the disturbance decoupling problem by static state feedback for a system ruled by $\mathcal{D}\,\x(t)=A\,\x(t)+B\,\u(t)+L\,\w(t)$ and $\y(t)=C\,\x(t)$.


\begin{theorem}
\label{the11}
Let $\ima L \subseteq \gV^\star_{\scriptscriptstyle (A,B,C,D)}$. The following hold:
\begin{description}
\item{{\em i)}} $\gV^\star_{\scriptscriptstyle (A,B,C,D)}={\gV}^\star_{\scriptscriptstyle (A,[\,B\;\;L\,],C,[\,D\;\; 0 \,])}$;
\item{{\em ii)}} $\Phi_{\scriptscriptstyle (A,[\,B\;\;L\,],C,[\,D\;\; 0 \,])}\subseteq \Phi_{\scriptscriptstyle (A,B,C,D)}$;
\item{{\em iii)}} $\forall\,{\gV}\in\Phi_{\scriptscriptstyle (A,[\,B\;\;L\,],C,[\,D\;\; 0 \,])}$,  $\quad\ima L   \subseteq \gV$.
\end{description}
\end{theorem}


\begin{theorem} 
\label{the12}
$\ima L \subseteq {\gV}^\star_{\scriptscriptstyle (A,B,C,D)}\;\; \Leftrightarrow \;\; \ima L \subseteq \gR^\star_{\scriptscriptstyle (A,[\,B\;\;L\,],C,[\,D\;\; 0 \,])}$.
\end{theorem}


\begin{theorem}
\label{the13}
If $\ima L \subseteq \gV^\star_{\scriptscriptstyle (A,B,C,D)}$,
the subspace $\gR^\star_{\scriptscriptstyle (A,[\,B\;\;L\,],C,[\,D\;\; 0 \,])}$ is the smallest of all the $(A,B,C,D)$-self bounded subspaces $\gV$ satisfying $\ima L \subseteq \gV$.
\end{theorem}


The following three results are a generalization of the last three: they are concerned with a geometric condition in the form $\ima \bsmat  L_1 \\[1mm] L_2 \esmat \subseteq (\gV^\star_{\scriptscriptstyle (A,B,C,D)} \oplus 0_{\scriptscriptstyle \scriptscriptstyle \gY})+\ima \bsmat B \\[1mm] D \esmat$.


\begin{theorem}
\label{the14}
Let $\ima \bsmat  L_1 \\[1mm] L_2 \esmat \subseteq (\gV^\star_{\scriptscriptstyle (A,B,C,D)} \oplus 0_{\scriptscriptstyle \scriptscriptstyle \gY})+\ima \bsmat B \\[1mm] D \esmat$.
The following results hold:
\begin{description}
\item{{\em i)}} $\gV^\star_{\scriptscriptstyle (A,B,C,D)}={\gV}^\star_{\scriptscriptstyle (A,[\,B\;\;L_1\,],C,[\,D\;\; L_2 \,])}$;
\item{{\em ii)}} $\Phi_{\scriptscriptstyle (A,[\,B\;\;L_1\,],C,[\,D\;\; L_2 \,])}\subseteq \Phi_{\scriptscriptstyle (A,B,C,D)}$;
\item{{\em iii)}} $\forall\,{\gV}\in\Phi_{\scriptscriptstyle (A,[\,B\;\;L_1\,],C,[\,D\;\; L_2 \,])}, \quad \ima \bsmat L_1 \\[1mm] L_2 \esmat \subseteq \gV \oplus 0_{\scriptscriptstyle \scriptscriptstyle \gY}+\ima \bsmat B \\[1mm] D \esmat$.
\end{description}
\end{theorem}


\begin{theorem}
$\ima \bsmat  L_1 \\[1mm] L_2 \esmat \subseteq (\gV^\star_{\scriptscriptstyle (A,B,C,D)} \oplus 0_{\scriptscriptstyle \scriptscriptstyle \gY})+\ima \bsmat B \\[1mm] D \esmat\;\; \Leftrightarrow \;\; \ima \bsmat  L_1 \\[1mm] L_2 \esmat \subseteq (\gR^\star_{\scriptscriptstyle (A,[\,B\;\;L_1\,],C,[\,D\;\;L_2\,])} \oplus 0_{\scriptscriptstyle \scriptscriptstyle \gY})+\ima \bsmat B \\[1mm] D \esmat$.
\end{theorem}

\begin{theorem}
\label{prop3}
If
$\ima \bsmat  L_1 \\[1mm] L_2 \esmat \subseteq (\gV^\star_{\scriptscriptstyle (A,B,C,D)} \oplus 0_{\scriptscriptstyle \scriptscriptstyle \gY})+\ima \bsmat B \\[1mm] D \esmat$,
the subspace $\gR^\star_{\scriptscriptstyle (A,[\,B\;\;L_1\,],C,[\,D\;\; L_2 \,])}$ is the smallest of all the $(A,B,C,D)$-self bounded subspaces $\gV$ satisfying
$\ima \bsmat L_1 \\[1mm] L_2 \esmat \subseteq (\gV \oplus 0_{\scriptscriptstyle \scriptscriptstyle \gY})+\ima \bsmat B\\[1mm] D \esmat$.
\end{theorem}

We now dualize all the previous results. The first three involve an inclusion in the form $\gS_{\scriptscriptstyle (A,B,C,D)}^\star\subseteq \ker M$.

\begin{theorem}
\label{the11d}
Let $\gS_{\scriptscriptstyle (A,B,C,D)}^\star\subseteq \ker M$. The following hold:
\begin{description}
\item{{\em i)}} $\gS^\star_{\scriptscriptstyle (A,B,C,D)}={\gS}^\star_{\scriptscriptstyle \left(A,B,\bsmat {\scriptscriptstyle C} \\ {\scriptscriptstyle M} \esmat,\bsmat {\scriptscriptstyle D} \\ {\scriptscriptstyle 0} \esmat\right)}$;
\item{{\em ii)}} $\Psi_{\scriptscriptstyle \left(A,B,\bsmat {\scriptscriptstyle C}\\ {\scriptscriptstyle M} \esmat,\bsmat {\scriptscriptstyle D} \\ {\scriptscriptstyle 0} \esmat\right)} \subseteq \Psi_{\scriptscriptstyle (A,B,C,D)}$;
\item{{\em iii)}} $\forall\,{\gS}\in\Psi_{\scriptscriptstyle \left(A,B,\bsmat {\scriptscriptstyle C}  \\{\scriptscriptstyle M} \esmat,\bsmat {\scriptscriptstyle D} \\ {\scriptscriptstyle 0} \esmat\right)}, \quad \gS   \subseteq \ker M$.
\end{description}
\end{theorem}


\begin{theorem}
\label{the12d}
${\gS}^\star_{\scriptscriptstyle (A,B,C,D)}\subseteq \ker M\; \Leftrightarrow\;\gQ^\star_{\scriptscriptstyle \left(A,B,\bsmat {\scriptscriptstyle C}\\ {\scriptscriptstyle M} \esmat,\bsmat {\scriptscriptstyle D} \\ {\scriptscriptstyle 0} \esmat\right)}\subseteq \ker M$.
\end{theorem}


\begin{theorem}
\label{the13d}
If $\gS^\star_{\scriptscriptstyle (A,B,C,D)} \subseteq \ker M$,
the subspace $\gQ^\star_{\scriptscriptstyle \left(A,B,\bsmat {\scriptscriptstyle C} \\ {\scriptscriptstyle M} \esmat,\bsmat {\scriptscriptstyle D} \\ {\scriptscriptstyle 0} \esmat\right)}$ is the largest of all the $(A,B,C,D)$-self hidden subspaces $\gS$ satisfying
$\gS\subseteq \ker M$.
\end{theorem}


Finally, we consider the generalization $(\gS_{\scriptscriptstyle (A,B,C,D)}^\star \oplus \gU) \cap \ker [\begin{array}{cc} C & D \end{array}]\subseteq \ker [\begin{array}{cc} M_1 & M_2 \end{array}]$ of the condition $\gS^\star_{\scriptscriptstyle (A,B,C,D)} \subseteq \ker M$.


\begin{theorem}
\label{the11d1}
Let $(\gS_{\scriptscriptstyle (A,B,C,D)}^\star \oplus \gU) \cap \ker [\begin{array}{cc} C & D \end{array}]\subseteq \ker [\begin{array}{cc} M_1 & M_2 \end{array}]$. 
 The following results hold:
\begin{description}
\item{{\em i)}} $\gS^\star_{\scriptscriptstyle (A,B,C,D)}={\gS}^\star_{\scriptscriptstyle \left(A,B,\bsmat {\scriptscriptstyle C} \\ {\scriptscriptstyle M_2} \esmat,\bsmat {\scriptscriptstyle D} \\ {\scriptscriptstyle M_2} \esmat\right)}$;
\item{{\em ii)}} $\Psi_{\scriptscriptstyle \left(A,B,\bsmat {\scriptscriptstyle C}\\{\scriptscriptstyle M_1} \esmat,\bsmat {\scriptscriptstyle D} \\ {\scriptscriptstyle M_2} \esmat\right)} \subseteq \Psi_{\scriptscriptstyle (A,B,C,D)}$;
\item{{\em iii)}} $\forall\,{\gS}\!\in\!\Psi_{\scriptscriptstyle \left(A,B,\bsmat {\scriptscriptstyle C} \\{\scriptscriptstyle M_1} \esmat,\bsmat {\scriptscriptstyle D} \\ {\scriptscriptstyle M_2} \esmat\right)}, \;(\gS \!  \oplus \!\gU) \cap \ker [\begin{array}{cc} C & D \end{array}]\!\subseteq \ker [\begin{array}{cc} M_1 & M_2 \end{array}]$.
\end{description}
\end{theorem}


\begin{theorem}
$(\gS_{\scriptscriptstyle (A,B,C,D)}^\star \oplus \gU) \cap \ker [\begin{array}{cc} C& D \end{array}]\subseteq \ker [\begin{array}{cc} M_1 & M_2\end{array}]\;\; \Leftrightarrow \;\; (\gQ^\star_{\scriptscriptstyle \left(A,B,\bsmat {\scriptscriptstyle C} \\ {\scriptscriptstyle M_1} \esmat,\bsmat {\scriptscriptstyle D} \\ {\scriptscriptstyle M_2} \esmat\right)} \oplus \gU) \cap \ker [\begin{array}{cc} C & D \end{array}]\subseteq \ker [\begin{array}{cc} M_1 & M_2 \end{array}]$.
\end{theorem}


\begin{theorem}
\label{the12d3}
If $(\gS_{\scriptscriptstyle (A,B,C,D)}^\star \oplus \gU) \cap \ker [\begin{array}{cc} C & D \end{array}]\subseteq \ker [\begin{array}{cc} M_1 & M_2 \end{array}]$, $\gQ^\star_{\scriptscriptstyle \left(A,B,\bsmat {\scriptscriptstyle C}\\ {\scriptscriptstyle M_1} \esmat,\bsmat {\scriptscriptstyle D} \\ {\scriptscriptstyle M_2} \esmat\right)}$ is the largest of all $(A,B,C,D)$-self hidden subspaces $\gS$ satisfying $(\gS \oplus \gU) \cap \ker [\begin{array}{cc}\! C\! & D\! \end{array}]\subseteq \ker [\begin{array}{cc} \!M_1\! & M_2\! \end{array}]$.
\end{theorem}

\end{document}